\documentclass[12pt,leqno]{article}

\usepackage{amsthm}
\usepackage[usenames]{color}

\usepackage{amsmath}
\usepackage{amscd}
\usepackage{amsmath}
\usepackage{amssymb}
\usepackage{latexsym}
\makeatletter

\@addtoreset{equation}{section}
\makeatother

\newtheorem{thm}{Theorem}[subsection]
\newtheorem{pr}[thm]{Proposition}
\newtheorem{df}[thm]{Definition}
\newtheorem{lm}[thm]{Lemma}
\newtheorem{cor}[thm]{Corollary}

\newtheorem{rmk}[thm]{Remark}

\newcommand{\sm}{\raisebox{2.33pt}{~\rule{6.4pt}{1.3pt}~}}

\input xy \xyoption{all} \CompileMatrices

\begin{document}

\title{Cotangent bundle
and micro-supports\\
in mixed characteristic case}
\author{Takeshi Saito}

\maketitle

\maketitle
\begin{abstract}
For a regular scheme
and a prime number $p$,
we define the FW-cotangent bundle
as a vector bundle
on the closed subscheme
defined by $p=0$,
under a certain finiteness condition.

For a constructible complex
on the \'etale site of the scheme,
we introduce the condition to be
micro-supported on a closed
conical subset in the FW-cotangent bundle.
At the end of the article,
we compute the singular supports
in some cases.
\end{abstract}

Let $k$ be a perfect field 
of characteristic $p>0$
and let $X$ be a regular 
noetherian scheme
such that the closed
subscheme
$X_{{\mathbf F}_p}$
defined by $p=0$
is a scheme of
finite type over $k$.
For example,
$X$ is of finite type
over a discrete
valuation ring 
${\cal O}_K$
with residue field $k$
or over $k$ itself.
The main purpose of the article is to prepare
a framework to study the
micro-support for an \'etale sheaf
on $X$
as in the transcendental setting \cite{KSc}
or in the setting of algebraic geometry
\cite{Be},
by introducing a variant of
the cotangent bundle defined on
$X_{{\mathbf F}_p}$
in an arithmetic setting.

The key property used in the definition
of singular supports in \cite{Be}
in the geometric case
is the local acyclicity of
morphisms 
to other smooth schemes
from the scheme
where the sheaf is defined.
However, as we see in Remark \ref{rmk},
a simple imitation of
the geometric case does not work
in the mixed characteristic situation
as we do not have sufficiently many
morphisms out of a scheme.
Instead, we will use 
the ${\cal F}$-transversality
introduced in \cite[Definition 8.5]{CC}
of morphisms from other
schemes to the scheme
where the sheaf is defined,
which is known to give a characterization
of singular support
in the geometric case.

Let $\Lambda$ be a finite field
of characteristic $\neq p$.
For a separated 
morphism $h\colon W\to X$
of finite type of regular noetherian
schemes 
and a constructible complex
${\cal F}$ of $\Lambda$-modules on 
the \'etale site of $X$,
we define the ${\cal F}$-transversality
(Definition \ref{dfFtr})
as the property 
for the canonical morphism
$c_{{\cal F},h}\colon
h^*{\cal F}\otimes
Rh^!\Lambda\to
Rh^!{\cal F}$
(\ref{eqcF})
to be an isomorphism,
similarly as in \cite[Definition 8.5]{CC}.
We show that the transversality
for the direct image
is equivalent to the property
that the base change morphism
is an isomorphism in
Proposition \ref{prhF}.2.

We have shown that
the sheaf $F\Omega^1_X$
of FW-differentials
is a locally free
${\cal O}_{X_{{\mathbf F}_p}}$-modules
of rank $\dim X$
in \cite{FW}.
We call the associated
vector bundle
$FT^*X|_{X_{{\mathbf F}_p}}$
on $X_{{\mathbf F}_p}$ the
FW-cotangent bundle of $X$.
The fiber 
$FT^*X|_x$
at a closed point
$x\in X_{{\mathbf F}_p}$
is canonically identified (\ref{eqTx})
with the Frobenius pull-back
$F^*T^*_xX$
of the cotangent space
that is the 
vector space
${\mathfrak m}_x/
{\mathfrak m}_x^2$
regarded as a scheme
over the residue field $k(x)$.
For a closed conical subset 
$C$ of the vector bundle
$FT^*X|_{X_{{\mathbf F}_p}}$
and for a morphism $h\colon W\to X$
of finite type of regular schemes,
we define the $C$-transversality
in Definition \ref{dfhC}.1
similarly as in \cite[1.2]{Be}.

Using the $C$-transversality
and the ${\cal F}$-transversality,
we define the condition for
${\cal F}$ to be micro-supported on $C$
in Definition \ref{dfms}. 
This is a property along the
closed subscheme 
$X_{{\mathbf F}_p}$.
For example,
if $X$ is of finite
type over ${\cal O}_K$ as above,
then 
${\cal F}$ is locally constant on 
a neighborhood of the closed fiber
$X_{{\mathbf F}_p}$ if and only if
${\cal F}$ is micro-supported on
the $0$-section 
$FT^*_XX|_{X_{{\mathbf F}_p}}$.

If the smallest closed conical
subset of $FT^*X|_{X_{{\mathbf F}_p}}$
on which ${\cal F}$ is micro-supported
exists,
we call it the singular
support $SS{\cal F}$ of ${\cal F}$.
The author does not know how to show
the existence in general.
We compute the singular
support of some 
sheaves on regular schemes
in Proposition \ref{prKum}.


The author thanks Luc Illusie for comments
on an earlier version.
The author thanks greatly
an anonymous referee
for careful reading and
for the suggestion of a generalization
of Lemma \ref{lmXW} and its proof 
as in the current form. 
The research is partially supported
by Grant-in-Aid (B) 19H01780.

%

\tableofcontents

\section{${\cal F}$-transversality}\label{sF}

In this section, we study properties
of morphisms of schemes
with respect to
complexes on the \'etale site
of a scheme.
The transversality is defined 
as a condition for a canonical morphism
for extraordinary pull-back to be
an isomorphism.
In Section \ref{ssFtr},
after preparing some sorites on
the canonical morphism,
we establish basic properties
on the transversality.
In Section \ref{ssla},
after recalling basic properties
of local acyclicity,
we study the relation between
the local acyclicity and
the transversality.

In this section
and Section \ref{sms},
$\Lambda$ denotes
a finite field of characteristic $\ell$
invertible on relevant 
noetherian schemes.
The derived categories
$D^+(-,\Lambda)$ of 
bounded below complexes
and 
$D^b_c(-,\Lambda)$ of 
constructible complexes
are defined as usual.

\subsection{${\cal F}$-transversality}\label{ssFtr}

Let $h\colon W\to X$
be a separated morphism
of finite type of noetherian schemes
and $\Lambda$ be a finite field
of characteristic $\ell$ invertible on $X$.
The functor
$Rh^!\colon D^+(X,\Lambda)
\to D^+(W,\Lambda)$
is defined as the adjoint
of
$Rh_!\colon D(W,\Lambda)
\to D(X,\Lambda)$
in \cite[Th\'eor\`eme 3.1.4.]{DP}.
If $X$ is quasi-excellent,
by the finiteness theorem
\cite[{\sc Th\'eor\`eme} 1.1.1]{fini},
we have a functor
$Rh^!\colon D^b_c(X,\Lambda)
\to D^b_c(W,\Lambda)$
see also \cite[Corollaire 1.5]{TF}.
Recall that a scheme of
finite type over
a Dedekind domain with
fraction field of characteristic 0
is quasi-excellent
by \cite[Scholie (7.8.3)]{EGA4}.

Let ${\cal F}
\in D^+(X,\Lambda)$ and 
${\cal G}\in D^+(W,\Lambda)$ .
Then, the adjoint of the morphism
$h^*{\cal F}\otimes h^*Rh_*{\cal G}
\to h^*{\cal F}\otimes {\cal G}$
induced by the adjunction
$h^*Rh_*{\cal G} \to{\cal G}$
defines a canonical morphism
\begin{equation}
{\cal F}\otimes Rh_*{\cal G}
\to Rh_*(h^*{\cal F}\otimes {\cal G}).
\label{eqpr0}
\end{equation}
If $h$ is an open immersion
and if ${\cal G}=h^*{\cal G}_X$
for some extension of ${\cal G}$
on $X$,
(\ref{eqpr0}) is identified with
the morphism
${\cal F}\otimes R{\cal H}om(h_!\Lambda,
{\cal G}_X)
\to R{\cal H}om(h_!\Lambda,{\cal F}\otimes {\cal G}_X)$ defined
by the product.

Applying the construction
(\ref{eqpr0})
to a compactification of $h$
and the extension by $0$,
a canonical isomorphism
\begin{equation}
{\cal F}\otimes
Rh_!{\cal G} 
\to Rh_!(h^*{\cal F}\otimes{\cal G} ),
\label{eqprj}
\end{equation}
the projection formula
\cite[(4.9.1)]{Rapport}
is defined.
\if{This is defined as the adjoint
$h^*{\cal F}\otimes^L_\Lambda
h^*Rh_!{\cal G} 
\to h^*{\cal F}\otimes^L_\Lambda{\cal G}$
of the morphism induced
by the adjunction
$h^*Rh_!{\cal G} 
\to {\cal G}$
if $h$ is proper.
It is defined as the inverse of
the isomorphism
$h^*{\cal F}\otimes^L_\Lambda
h^*Rh_!{\cal G} 
\gets h^*{\cal F}\otimes^L_\Lambda{\cal G}$
if $h$ is an open immersion.}\fi

\begin{df}\label{dfAB}
Let $h\colon W\to X$
be a separated morphism of finite
type of quasi-excellent noetherian schemes.
Let ${\cal F}\in D^+(X,\Lambda)$.

{\rm 1.}
Let ${\cal G}\in D^+(X,\Lambda)$. 
We define a canonical morphism
\begin{equation}
c_{{\cal F},{\cal G},h}\colon
h^*{\cal F}
\otimes
Rh^!{\cal G}
\to 
Rh^!({\cal F}
\otimes
{\cal G})
\label{eqAB}
\end{equation}
to be the adjoint of the composition
$$
Rh_!(h^*{\cal F}
\otimes
Rh^!{\cal G})
\to 
{\cal F}
\otimes
Rh_!Rh^!{\cal G}
\to {\cal F}
\otimes
{\cal G}$$
of the inverse 
of the isomorphism {\rm (\ref{eqprj})}
and the morphism induced
by the adjunction
$Rh_!Rh^!{\cal G}
\to
{\cal G}$.
For ${\cal G}=\Lambda$,
we define a canonical morphism
\begin{equation}
c_{{\cal F},h}
\colon 
h^*{\cal F}
\otimes^L
Rh^!\Lambda
\to Rh^!{\cal F}
\label{eqcF}
\end{equation}
to be 
$c_{{\cal F},\Lambda,h}$.
\end{df}

\begin{lm}\label{lmcF}
Let $h\colon W\to X$
be a separated morphism of finite
type of noetherian schemes.
Let ${\cal F}\in D^+(X,\Lambda)$.

{\rm 1.}
Let ${\cal G},{\cal H}\in D^+(X,\Lambda)$.
Then, the diagram 
\begin{equation}
\begin{CD}
h^*{\cal F}
\otimes 
Rh^!({\cal G}\otimes {\cal H})
@>{c_{{\cal F},
{\cal G}\otimes {\cal H},h}}>>
Rh^!({\cal F}
\otimes 
{\cal G}\otimes {\cal H})\\
@A{1\otimes c_{{\cal G},
{\cal H},h}}AA
@AA{c_{{\cal F},
{\cal G},h}\otimes 1}A\\
h^*{\cal F}
\otimes 
Rh^!{\cal G}\otimes h^*{\cal H}
@>{c_{{\cal F},
{\cal G},h}\otimes 1}>>
Rh^!({\cal F}\otimes {\cal G})
\otimes 
h^*{\cal H}
\end{CD}
\end{equation}
is commutative.

{\rm 2.}
Let $g\colon V\to W$
be a separated morphism of finite
type of schemes
and 
let ${\cal G}\in D^+(X,\Lambda)$.
Then, the diagram
\begin{equation}
\xymatrix{
(hg)^*{\cal F}
\otimes
R(hg)^!{\cal G}
\ar[r]^{c_{{\cal F},{\cal G},hg}}
&
R(hg)^!({\cal F}\otimes
{\cal G})
\\
g^*h^*{\cal F}
\otimes
Rg^!Rh^!{\cal G}
\ar[u]
\ar[rd]^{c_{h^*{\cal F},
Rh^!{\cal G},g}}
&
Rg^!Rh^!({\cal F}\otimes
{\cal G})
\ar[u]
\\
g^*h^*{\cal F}
\otimes
Rg^!\Lambda
\otimes
g^*Rh^!{\cal G}
\ar[u]^{1\otimes
c_{Rh^!{\cal G},g}}
\ar[dr]_{
c_{h^*{\cal F},g}\otimes 1}
&
Rg^!(h^*{\cal F}
\otimes Rh^!{\cal G})
\ar[u]_{Rg^!(c_{{\cal F},{\cal G},h})}
\\
&
Rg^!h^*{\cal F}
\otimes g^*Rh^!{\cal G}.
\ar[u]_{c_{h^*{\cal F},Rh^!{\cal G},g}}
}
\label{eqcgh}
\end{equation}
where the upper vertical arrows
are canonical isomorphisms
{\rm \cite[(3.1.13.1)]{DP}}
is commutative.

{\rm 3.}
Let $$\begin{CD}
X@<h<< W\\
@VfVV@VV{f'}V\\
Y@<g<<V
\end{CD}$$
be a cartesian diagram of
separated morphisms
of finite type.
Then, the diagram
\begin{equation}
\begin{CD}
g^*Rf_*{\cal F}
\otimes
Rg^!\Lambda
@>{c_{Rf_*{\cal F},g}}>>
Rg^!Rf_*{\cal F}
\\
@VVV@VVV\\
Rf'_*h^*{\cal F}
\otimes
Rg^!\Lambda
@.
Rf'_*Rh^!{\cal F}
\\
@V
{\rm (\ref{eqpr0})}VV
@AA{Rf'_*(c_{{\cal F},h})}A\\
Rf'_*(h^*{\cal F}
\otimes
f'^*Rg^!\Lambda)
@>>>
Rf'_*(h^*{\cal F}
\otimes
Rh^!\Lambda)
\end{CD}
\label{eqcfg}
\end{equation}
where the arrows without
tags are defined by
base change morphisms
is commutative.
\end{lm}

\proof{
1.
The diagram
$$\begin{CD}
Rh_!Rh^!({\cal G}\otimes {\cal H})
@>>>
{\cal G}\otimes {\cal H}\\
@A
{Rh_!(c_{{\cal G},{\cal H},h})}AA
@AAA\\
Rh_!(Rh^!{\cal G}\otimes h^*{\cal H})
@<{{\rm (\ref{eqprj})}}<<
Rh_!Rh^!{\cal G}\otimes {\cal H}
\end{CD}$$
where the arrows without
tags are defined by the adjunction
is commutative
by the definition of
$c_{{\cal G},{\cal H},h}$.
Taking the tensor products with ${\cal F}$,
applying the projection formula
(\ref{eqprj}) and
taking the adjoint,
we see that the upper triangle in
\begin{equation*}
\xymatrix{
h^*{\cal F}
\otimes 
Rh^!({\cal G}\otimes {\cal H})
\ar[r]^{c_{{\cal F},
{\cal G}\otimes {\cal H},h}}&
Rh^!({\cal F}
\otimes 
{\cal G}\otimes {\cal H})\\
h^*{\cal F}
\otimes 
Rh^!{\cal G}\otimes h^*{\cal H}
\ar[u]^{1\otimes c_{{\cal G},
{\cal H},h}}
\ar[ru]^
{c_{{\cal F}\otimes{\cal H},
{\cal G},h}}
\ar[r]^{c_{{\cal F},
{\cal G},h}\otimes 1}
&
Rh^!({\cal F}\otimes {\cal G})
\otimes 
h^*{\cal H}
\ar[u]_{c_{{\cal F},
{\cal G},h}\otimes 1}
}
\end{equation*}
is commutative.
The lower triangle is similarly
commutative
and the assertion follows.

2.
The lower quadrangle
is commutative by 1.
The composition 
$g^*h^*{\cal F}
\otimes
Rg^!Rh^!{\cal G}
\to
Rg^!Rh^!{\cal F}$
through
$
Rg^!(h^*{\cal F}
\otimes Rh^!{\cal G})$
is the adjoint of
$Rh_!Rg_!
(g^*h^*{\cal F}
\otimes
Rg^!Rh^!{\cal G})
\to
{\cal F}\otimes
Rh_!Rg_!
Rg^!Rh^!{\cal G}$
induced by
the adjunction
$Rh_!Rg_!
Rg^!Rh^!{\cal G}
\to
Rh_!Rh^!{\cal G}
\to {\cal G}$.
Since the last morphism
is identified
with 
the adjunction
$R(hg)_!
R(hg)^!{\cal G}
\to {\cal G}$,
the upper pentagon is also commutative.

3.
For ${\cal G}\in D^+(V,\Lambda)$,
we consider the diagram
\begin{equation}
\begin{CD}
f^*Rg_!(g^*Rf_*{\cal F}
\otimes
{\cal G})
@<{f^*{\rm (\ref{eqprj})}}<<
f^*Rf_*{\cal F}
\otimes
f^*Rg_!{\cal G}
@>>>
{\cal F}
\otimes
f^*Rg_!{\cal G}
\\
@VVV@.@VVV\\
Rh_!f'^*(Rf'_*h^*{\cal F}
\otimes
{\cal G})
@>>>
Rh_!(h^*{\cal F}
\otimes
f^*{\cal G})
@<{\rm (\ref{eqprj})}<<
{\cal F}
\otimes
Rh_!f'^*{\cal G}
\end{CD}
\label{eqcfga}
\end{equation}
defined as follows.
The vertical arrows are
defined by the base change morphisms
and the horizontal arrows
without labels are
defined by adjunction.
We see that the diagram is commutative
by reducing to the case
where $g$ is proper and
going back to the definition
of (\ref{eqprj}).

We apply (\ref{eqcfga}) to
${\cal G}=Rg^!\Lambda$.
Since the composition
$f^*Rg_!Rg^!\Lambda
\to
Rh_!f'^*Rg^!\Lambda
\to Rh_!Rh^!\Lambda
\to \Lambda$
of the base change morphisms
with the adjunction
is induced by the adjuncion
$Rg_!Rg^!\Lambda\to \Lambda$,
we obtain a commutative diagram
\begin{equation}
\begin{CD}
f^*Rg_!(g^*Rf_*{\cal F}
\otimes
Rg^!\Lambda)
@<{f^*{\rm (\ref{eqprj})}}<<
f^*Rf_*{\cal F}
\otimes
f^*Rg_!Rg^!\Lambda
@>>>
{\cal F}
\\
@VVV@.@AAA\\
Rh_!f'^*(Rf'_*h^*{\cal F}
\otimes
Rg^!\Lambda)
@>>>
Rh_!(h^*{\cal F}
\otimes
Rh^!\Lambda)
@<{\rm (\ref{eqprj})}<<
{\cal F}
\otimes
Rh_!Rh^!\Lambda
\end{CD}
\label{eqcfgb}
\end{equation}
Since the canonical morphism
(\ref{eqcF}) is defined as
the adjoint of (\ref{eqprj}),
we obtain (\ref{eqcfg})
by taking the adjoint of (\ref{eqcfgb}).
\qed

}

\begin{lm}\label{lmij}
Let $i\colon Z\to X$ be a closed
immersion of noetherian schemes
and let ${\cal F},{\cal G}
\in D^+(X,\Lambda)$.

{\rm 1.}
We define the slant arrow
and the vertical arrow
in the diagram
\begin{equation}
\xymatrix{
{\cal F}\otimes
i_*Ri^!{\cal G}
\ar[r]^-{\rm(\ref{eqprj})}
\ar[rd]
&
i_*(i^*{\cal F}\otimes
Ri^!{\cal G})
\ar[r]^-{i_*(c_{{\cal F},{\cal G},i})}
&
i_*Ri^!({\cal F}\otimes {\cal G})
\ar[d]
\\
&
{\cal F}\otimes
R{\cal H}om(i_*\Lambda,{\cal G})
\ar[r]
&
R{\cal H}om(i_*\Lambda,{\cal F}
\otimes{\cal G})
}
\label{eqic}
\end{equation}
by the canonical isomorphism
$i_*Ri^!\to R{\cal H}om(i_*\Lambda,-)$
and the lower horizontal arrow
by the product.
Then, the diagram
{\rm (\ref{eqic})}
is commutative.

{\rm 2.}
Let $j\colon U=X\sm Z\to X$
be the open immersion of the complement.
Then, 
the exact sequence
$0\to j_!\Lambda
\to \Lambda\to i_*\Lambda\to 0$
defines a commutative diagram
\begin{equation}
\begin{CD}
{\cal F}\otimes
i_*Ri^!{\cal G}
@>>> {\cal F}\otimes
{\cal G}
@>>>
{\cal F}\otimes
Rj_*j^*{\cal G}
@>>>\\
@V{c_{{\cal F},{\cal G},i}}VV@|
@VV{\rm (\ref{eqpr0})}V@.\\
i_*Ri^!({\cal F}\otimes
{\cal G})
@>>> {\cal F}\otimes
{\cal G}
@>>>
Rj_*j^*({\cal F}\otimes
{\cal G})
@>>>
\end{CD}
\label{eqij}
\end{equation}
of distinguished triangles.
\end{lm}

\proof{
1.
By the definition of
$c_{{\cal F},{\cal G},i}$,
the morphism
$i_*(c_{{\cal F},{\cal G},i})\colon
i_*(i^*{\cal F}\otimes
Ri^!{\cal G})
\to i_*Ri^!({\cal F}
\otimes {\cal G})$
is the unique morphism
such that the diagram
$$
\begin{CD}
{\cal F}\otimes
i^*Ri^!{\cal G}
@>>> {\cal F}\otimes{\cal G}
\\
@V{\rm (\ref{eqprj})}VV@AAA\\
i_*(i^*{\cal F}\otimes
Ri^!{\cal G})
@>{i_*(c_{{\cal F},{\cal G},i})}>>
i_*Ri^!({\cal F}\otimes{\cal G})
\end{CD}$$
is commutative.
Here the arrows without tag
are defined by the
adjunction $i_*Ri^!\to 1$.
Similarly, the lower horizontal arrow
${\cal F}\otimes
R{\cal H}om(i_*\Lambda,{\cal G})
\to
R{\cal H}om(i_*\Lambda,{\cal F}
\otimes{\cal G})$
is the unique morphism
such that the diagram
$$
\begin{CD}
{\cal F}\otimes
i^*Ri^!{\cal G}
@>>> {\cal F}\otimes{\cal G}
\\
@VVV@AAA\\
{\cal F}\otimes
R{\cal H}om(i_*\Lambda,{\cal G})
@>>>
R{\cal H}om(i_*\Lambda,{\cal F}
\otimes{\cal G})
\end{CD}$$
is commutative.
Here the left vertical arrow
is the slant arrow in (\ref{eqic})
and the right vertical arrow
is induced by $\Lambda\to i_*\Lambda$.
Hence
the assertion follows.

2.
The exact sequence
$0\to j_!\Lambda
\to \Lambda\to i_*\Lambda\to 0$
defines a commutative diagram
\begin{equation}
\begin{CD}
{\cal F}\otimes
R{\cal H}om(i_*\Lambda,{\cal G})
@>>> {\cal F}\otimes
{\cal G}
@>>>
{\cal F}\otimes
R{\cal H}om(j_!\Lambda,{\cal G})
@>>>\\
@VVV@VVV@VVV@.\\
R{\cal H}om(i_*\Lambda,{\cal F}\otimes
{\cal G})
@>>> {\cal F}\otimes
{\cal G}
@>>>
R{\cal H}om(j_!\Lambda,{\cal F}\otimes
{\cal G})
@>>>
\end{CD}
\label{eqij2}
\end{equation}
of distinguished triangles.
By 1.,
the left vertical arrow 
of (\ref{eqij}) is
identified with
that of (\ref{eqij2})
and similarly for the
right vertical arrows.
\qed

}

\begin{lm}\label{lmtrbc}
Let
$$\begin{CD}
X@<h<<W\\
@VfVV@VV{f'}V\\
Y@<g<<V
\end{CD}$$
be a cartesian of 
morphisms of finite type
of regular noetherian schemes.
If $f$ and $g$ are transversal,
then the base change
morphism
$f^*Rg^!\Lambda\to
Rh^!\Lambda$
is an isomorphism
of locally constant complexes.
\end{lm}

We say that a complex
${\cal F}$ is locally constant,
if its cohomology sheaf
${\cal H}^q{\cal F}$ is locally constant
for every $q$
and 
if ${\cal H}^q{\cal F}=0$
except for finitely many $q$.

\proof{
Since the assertion is local,
we may assume that
the morphism $g$ is a composition
of a smooth morphism
and a regular immersion.
Hence, it suffices to show
each case.

Assume that $g$ is a smooth
of relative dimension $d$.
Then, the adjoint
of the trace morphism
$Rg_!\Lambda(d)[2d]\to
\Lambda$ 
\cite[Th\'eor\`eme 2.9]{DP} defines
an isomorphism
$\Lambda(d)[2d]\to
Rg^!\Lambda$ by Poincar\'e duality
\cite[Th\'eor\`eme 3.2.5]{DP}.
Since the formation of
the trace morphism commutes
with base change,
the assertion follows in this case.

Assume that $g$ is a regular immersion 
of codimension $c$.
Then, by the absolute purity \cite[{\sc Th\'eor\`eme 3.1.1}]{purete},
the fundamental
class $[V]$ defines an isomorphism
$\Lambda\to Rg^!
\Lambda(c)[2c]$.
Since $f$ and $g$ 
are transversal,
further by the absolute purity,
the fundamental
class $[W]=f^![V]$
defines an isomorphism
$\Lambda\to Rh^!
\Lambda(c)[2c]$.
Hence the base change morphism
$f^*Rg^!
\Lambda
\to Rh^!\Lambda$ is an isomorphism.
\qed

}

\begin{df}\label{dfFtr}
Let $h\colon W\to X$
be a separated morphism of finite
type of noetherian schemes
and let ${\cal F}\in D^+(X,\Lambda)$.
We say that $h$ is ${\cal F}$-transversal
if the canonical morphism
{\rm (\ref{eqcF})} is an isomorphism.
\end{df}

For a closed immersion
$i\colon Z\to X$ of
regular noetherian schemes
and a separated morphism
$h\colon W\to X$
of finite type of
regular noetherian schemes,
we show that
$h$ is $i_*\Lambda$-transversal
if $h$ and $i$ are transversal
in Corollary \ref{cortrC}.2.
If $h$ is also an immersion,
if $Z$ and $W$ meets properly,
if the reduced part
$V$ of $Z\times_XW$ is regular
and if
the intersection multiplicity
$\mu(Z,W)$ 
is invertible in $\Lambda$,
then 
$h$ is still $i_*\Lambda$-transversal.
Hence the converse does not hold.

\begin{lm}\label{lmPoi}
Let $h\colon W\to X$
be a separated morphism of finite
type of noetherian schemes
and let ${\cal F}\in D^+(X,\Lambda)$.

{\rm 1.}
If $h\colon W\to X$
is smooth,
then $h$ is ${\cal F}$-transversal.

{\rm 2.}
If ${\cal F}$ is locally constant,
then $h$ is ${\cal F}$-transversal.
\end{lm}

\proof{
1.
This is exactly the Poincar\'e duality
\cite[Th\'eor\`eme 3.2.5]{DP}.

2.
Since the assertion is \'etale local,
the assertion is reduced to the case
where ${\cal F}=\Lambda$
by devissage.
\qed

}

\begin{lm}\label{lmiZ}
Let $i\colon Z\to X$ be a closed
immersion of noetherian schemes
and let ${\cal F}\in D^+(X,\Lambda)$.

{\rm 1.}
Assume that
$Z$ is the union of
closed subsets $Z_1,\ldots,Z_n
\subset X$
and that
for each subset $I\subset
\{1,\ldots,n\}$,
the immersion
$i_I\colon Z_I=\bigcap_{j\in I}Z_j\to X$
is ${\cal F}$-transversal.
Then, $i\colon Z\to X$
is ${\cal F}$-transversal.

{\rm 2.}
Let $j\colon U=X\sm Z\to X$
be the open immersion of the complement.
Then, the following conditions
are equivalent:

{\rm (1)}
$i\colon Z\to X$
is ${\cal F}$-transversal.

{\rm (2)}
The canonical morphism
${\cal F}
\otimes Rj_*\Lambda
\to
Rj_*j^*{\cal F}$
is an isomorphism.
\end{lm}

\proof{
1.
The quasi-isomorphism
$i_*\Lambda
\to[
\bigoplus_{j\in I}i_{j*}\Lambda\to
 \cdots \to
\bigoplus_{|I|=p}i_{I*}\Lambda
\to \cdots]$
defines a spectral sequence
$E_1^{p,q}=\bigoplus_{|I|=-p}
R^q{\cal H}om
(i_{I*}\Lambda,-)
\Rightarrow
R^{p+q}{\cal H}om
(i_*\Lambda,-)$.
By Lemma \ref{lmij}.1,
the assumption
implies that the morphisms
${\cal F}
\otimes R^q{\cal H}om
(i_{I*}\Lambda,\Lambda)
\to
R^q{\cal H}om
(i_{I*}\Lambda,{\cal F})$
on $E_1$-terms are isomorphisms.
Hence the assertion follows.

2.
The assertion follows from
Lemma \ref{lmij}.2 for ${\cal G}=\Lambda$.
\qed

}

\begin{pr}\label{prhF}
Let $h\colon W\to X$
be a separated morphism of finite
type of noetherian schemes
and let ${\cal F}\in D^+(X,\Lambda)$.
Assume that $h$ is ${\cal F}$-transversal.

{\rm 1}.
Assume that $Rh^!\Lambda$ 
is locally constant on
a neighborhood $W_1$
of ${\rm supp}\, h^*{\cal F}$
such that
$W_1\subset
{\rm supp}\, Rh^!\Lambda$.
Then, for a separated morphism 
$g\colon V\to W$
of finite type of noetherian schemes,
the following conditions are equivalent:

{\rm (1)}
$g$ is $h^*{\cal F}$-transversal.

{\rm (2)}
$hg$ is ${\cal F}$-transversal.

{\rm 2}.
Let 
$$\begin{CD}
X@<h<< W\\
@VfVV@VV{f'}V\\
Y@<g<<V
\end{CD}$$
be a cartesian diagram
of morphisms of
finite type of noetherian schemes.
Assume that $g$ is separated
and that
$Rg^!\Lambda$ is locally constant
of support $V$.
Further assume that 
the base change
morphism
\begin{equation}
f'^*Rg^!\Lambda\to
Rh^!\Lambda
\label{eqbcg}
\end{equation} 
is an isomorphism on 
a neighborhood of 
${\rm supp}\, h^*{\cal F}$.
Then the following conditions
are equivalent:

{\rm (1)}
The morphism
$g\colon V\to Y$
is $Rf_*{\cal F}$-transversal.

{\rm (2)}
The base change morphism
\begin{equation}
g^*Rf_*{\cal F}
\to Rf'_*h^*{\cal F}
\label{eqbch}
\end{equation}
is an isomorphism.
\end{pr}

By Lemma \ref{lmPoi}.1,
Proposition \ref{prhF}.2
(1)$\Rightarrow$(2)
gives a generalization of
the smooth base change theorem
\cite[Corollaire 1.2]{smbc}.

\proof{
1.
We consider the commutative diagram
(\ref{eqcgh}) for ${\cal G}=\Lambda$.
Since $h$ is assumed to
be ${\cal F}$-transversal,
$c_{{\cal F},h}$ is an isomorphism.
Since $Rh^!\Lambda$ is 
locally constant on a neighborhood
of the support of $h^*{\cal F}$,
the morphisms
$1\otimes c_{Rh^!\Lambda,g}$
and
$c_{h^*{\cal F},Rh^!\Lambda,g}$
are isomorphisms.
Hence
$c_{{\cal F},hg}$ is an isomorphism
if and only if
$c_{h^*{\cal F},g}\otimes 1$
is an isomorphisms.
Further by the assumption on
the support of $Rh^!\Lambda$,
the latter condition is equivalent
to the condition that
$c_{h^*{\cal F},g}$
is an isomorphism.

2.
We consider the commutative diagram
(\ref{eqcfg}).
By the proper base change theorem
or \cite[Corollaire 3.1.12.3]{DP},
the upper right vertical arrow
is an isomorphism.
Since $h$ is assumed ${\cal F}$-transversal,
the lower right vertical arrow
$Rf'_*(c_{{\cal F},h})$
is an isomorphism.
By the assumption on
$Rg^!\Lambda$,
the upper left vertical arrow
is an isomorphism if and only if
(\ref{eqbch}) is an isomorphism.
Further the arrow labeled
(\ref{eqpr0}) is an isomorphism.
Since (\ref{eqbcg}) 
is assumed to be an isomorphism,
the bottom horizontal arrow is an
isomorphism.
Hence the assertion follows
from the commutative diagram
(\ref{eqcfg}).
\qed

}

\begin{cor}\label{cortrC}
{\rm 1.}
Let the assumption
be the same as in Proposition
{\rm \ref{prhF}.2.}
Assume further that $f$ is proper
on the support of ${\cal F}$.
Then, $g$ is $Rf_*{\cal F}$-transversal.

{\rm 2.}
Let $p\colon Z\to X$ be 
a proper morphisms of regular schemes
and
let $h\colon W\to X$
be a separated morphism
of finite type of regular schemes.
If $h$ and $p$ are transversal,
then $h$ is $Rp_*\Lambda$-transversal.
\end{cor}

\proof{
1.
By the assumption that $f$ is proper
on the support of ${\cal F}$,
the base change morphism
(\ref{eqbch}) is an isomorphism
by the proper base change theorem.
Hence the assertion follows
from Proposition \ref{prhF}.2
(2)$\Rightarrow$(1).

2.
Let 
$$\begin{CD}
Z@<g<< V\\
@VpVV@VV{p'}V\\
X@<h<< W
\end{CD}$$
be a cartesian diagram.
By Lemma \ref{lmtrbc},
the base change morphism
$p'^*Rh^!
\Lambda
\to Rg^!\Lambda$ is an isomorphism.
Since $X$ and $W$ are regular,
the assumption in Proposition  \ref{prhF}.2
that
$Rh^!\Lambda$ is locally constant 
of support $W$ is satisfied
by the absolute purity \cite[{\sc Th\'eor\`eme 3.1.1}]{purete}
as in the proof of Lemma \ref{lmtrbc}.
Since $g$ is $\Lambda$-transversal
and $p$ is proper,
$h$ is $Rp_*\Lambda$-transversal
by 1.
\qed

}

\subsection{Local acyclicity
and ${\cal F}$-transversality}\label{ssla}

Let $f\colon X\to Y$
be a morphism of schemes
and $x$ and $y$ be
geometric points of $X$ and $Y$.
Let $f(x)$ denote the geometric point
of $Y$ defined by the composition
$x\to X\to Y$
and let
$X_{(x)}\to Y_{(f(x))}$ be the induced
morphism of strict localizations.
We call a morphism $y\to Y_{(f(x))}$
of schemes over $Y$
a specialization $f(x)\gets y$
and call $X_{(x),y}
=X_{(x)}
\times_{Y_{(f(x))}}y$
the Milnor fiber.
For a complex ${\cal F}$ of 
$\Lambda$-modules
on $X$,
the pull-back by $X_{(x),y}
\to X_{(x)}$
defines a canonical morphism
\begin{equation}
{\cal F}_x
=R\Gamma(X_{(x)},{\cal F}|_{{X_{(x)}}})
\to
R\Gamma(X_{(x),y},{\cal F}|_{{X_{(x),y}}}).
\label{eqMil}
\end{equation}

\begin{df}[{cf.~\cite[D\'efinition 2.12]{TF}}]\label{dfla}
Let $f\colon X\to Y$
be a morphism of schemes
and $Z\subset X$ be a closed
subset.
Let ${\cal F}$
be a complex of $\Lambda$-modules 
on $X$.
We say that $f$ is locally acyclic
relatively to ${\cal F}$
or ${\cal F}$-acyclic 
for short along $Z$ 
if for every geometric point
$x$ of $Z$ and for every specialization
$f(x)\gets y$,
the canonical morphism
${\cal F}_x\to R\Gamma(X_{(x),y},
{\cal F}|_{X_{(x),y}})$ 
{\rm (\ref{eqMil})} is an isomorphism.
If $X=Z$, we drop along $Z$
in the terminology.

We say that $f$ is universally
${\cal F}$-acyclic along $Z$, if
for every morphism $Y'\to Y$,
the base change
$X'\to Y'$ is locally acyclic 
relatively to the pull-back of
${\cal F}$ along the inverse image
$Z'\subset X'$ of $Z$.
\end{df}

\begin{lm}\label{lmlac}
Let $f\colon X\to Y$
be a morphism of schemes
and $Z\subset X$ be a closed
subset.
Let ${\cal F}
\in D^+(X,\Lambda)$.

{\rm 1.}
The following conditions are
equivalent.

{\rm (1)}
$f$ is ${\cal F}$-acyclic along $Z$.

{\rm (2)}
Let $s\gets t$ be a specialization
of geometric points of $Y$
such that $t$ is
the spectrum of an algebraic
closure of the residue field
of the point of $Y$ below $t$.
Let $Y_{(s)}$ denote the strict localization
and let
\begin{equation}
\begin{CD}
X_s@>{i'_s}>>X\times_YY_{(s)}@<{j'_t}<< X_t\\
@V{f_s}VV@V{f_{(s)}}VV @VV{f_t}V\\
s@>{i_s}>>Y_{(s)}@<{j_t}<<t
\end{CD}
\label{eqXYtpi}
\end{equation}
be the cartesian diagram.
Then, the canonical morphism 
\begin{equation}
i^{\prime*}_s{\cal F}
\to 
i^{\prime*}_sRj'_{t*}j^{\prime*}_t{\cal F}
\label{eqijst}
\end{equation}
is an isomorphism on 
the inverse image of $Z$.

{\rm 2.}
For a proper morphism $p\colon X\to P$
of schemes over $Y$,
we consider the following conditions:

{\rm (1)}
$f$ is ${\cal F}$-acyclic along $Z$.

{\rm (2)}
The morphism $g\colon P\to Y$ is 
$Rp_*{\cal F}$-acyclic along $p(Z)$.

\noindent
We have {\rm (1)}$\Rightarrow${\rm (2)}
 if $Z=p^{-1}(p(Z))$.
If $p$ is finite, we have
{\rm (2)}$\Rightarrow${\rm (1)}.

{\rm 3.}
The following conditions are
equivalent:

{\rm (1)}
$f$ is universally
${\cal F}$-acyclic along $Z$.

{\rm (2)}
For every 
smooth morphism $Y'\to Y$
and 
for the pull-back ${\cal F}'$
of ${\cal F}$ on 
$X'=X\times_YY'$,
the base change
$f'\colon X'\to Y'$ is 
${\cal F}'$-acyclic along 
the inverse image $Z'\subset X'$ of $Z$.
\end{lm}

Since the local acyclicity is a local property,
by locally taking an immersion
$X\to {\mathbf A}^n_Y$,
the study of local acyclicity
is reduced to the case
where $f\colon X\to Y$
is the projection
${\mathbf A}^n_Y\to Y$
by Lemma \ref{lmlac}.2.

\proof{
1.
A morphism $y'\to y$ of 
geometric points of $Y$
is the composition of a limit of 
smooth morphisms and
a homeomorphism in \'etale topology.
Hence for a geometric point $x$ of
$X$ and a specialization $f(x)\gets y$,
the pull-back
$R\Gamma(X_{(x),y},
{\cal F}|_{X_{(x),y}})\to
R\Gamma(X_{(x),y'},
{\cal F}|_{X_{(x),y'}})$
is an isomorphism by the smooth base
change theorem
\cite[Corollaire 1.2]{smbc}.
Thus, in the definition of
local acyclicity,
it suffices to consider 
specializations $f(x)\gets y$
such that $y$ is
the spectrum of an algebraic
closure of the residue field
of the point of $Y$ below $y$.

In the notation of (2),
for a geometric point $x$ of $X_s$,
the morphism induced by
(\ref{eqijst}) on the stalks 
at $x$ equals 
${\cal F}_x\to
R\Gamma(X_{(x),t},
{\cal F}|_{X_{(x),t}})$
(\ref{eqMil}).
Hence, the assertion follows.

2.
{\rm (1)}$\Rightarrow${\rm (2)}:
Let
\begin{equation}
\begin{CD}
X_s@>{i'_s}>>X\times_YY_{(s)}@<{j'_t}<< X_t\\
@V{p_s}VV@V{p_{(s)}}VV @VV{p_t}V\\
P_s@>{i''_s}>>P\times_YY_{(s)}@<{j''_t}<<P_t
\end{CD}
\end{equation}
be the base change of $X\to P$.
Then, the isomorphism (\ref{eqijst})
on the inverse image of $Z=p^{-1}(p(Z))$
implies an isomorphism
$i^{\prime\prime*}_sRp_{(s)*}{\cal F}
\to 
i^{\prime\prime*}_sRj''_{t*}j^{\prime\prime*}_tRp_{(s)*}{\cal F}$
on the inverse image of $p(Z)$
by proper base change theorem.

{\rm (2)}$\Rightarrow${\rm (1)}:
Let $z$ be a geometric point
of $Z$ and let $w=p(z)\gets y$
be a specialization.
Then the cospecialization morphism
$p_*{\cal F}_w
\to R\Gamma(P_{(w),y},p_*{\cal F}|_{P_{(w),y}})$
is the direct sum of
(\ref{eqMil}) for
$x\in p^{-1}(w)$
since $p$ is finite.
Hence
the assertion follows.

3.
Since the local acyclicity
is a local property
preserved by base change
by immersions
and commutes with limits,
the assertion follows.
\qed

}

\begin{lm}\label{lmla}
Let $X$ be a noetherian scheme
and ${\cal F}\in D^b(X,\Lambda)$.

{\rm 1.}
If ${\cal F}$ is locally constant
and if $f\colon X\to Y$
is smooth,
then $f$ is ${\cal F}$-acyclic.

{\rm 2.}
If $1_X\colon X\to X$
is ${\cal F}$-acyclic along $Z$
and if ${\cal F}$ is constructible,
then ${\cal F}$ is locally constant
on a neighborhood of $Z$.
\end{lm}

\proof{1.
By devissage,
the assertion follows from
the local acyclicity of smooth morphism
\cite[Th\'eor\`eme 2.1]{Artin}.

2.
For every geometric point  $s$ of
$Z$ and
every specialization 
$s\gets t$ of geometric points
of $X$,
the cospecialization morphism
${\cal F}_s\to {\cal F}_t$
is an isomorphism.
Hence the constructible sheaf
${\cal H}^q{\cal F}$ is locally constant
on a neighborhood of 
 every geometric point $s$ 
of $Z$ for every $q\in {\mathbf Z}$
by \cite[Proposition 2.11]{cst}.
Hence ${\cal F}$
is locally constant on a
neighborhood of $Z$.
\qed

}

\begin{pr}[{cf.~\cite[Proposition 8.11]{CC}}]\label{prla1}
Let $f\colon X\to Y$
be a smooth morphism
of regular schemes
of finite type over 
a discrete valuation ring
${\cal O}_K$
and $Z\subset X$
be a closed subset.
Let
${\cal F}$ be a constructible
complex of $\Lambda$-modules
on $X$.
Assume that 
for every separated morphism
$V\to Y$ of
regular schemes of finite
type over ${\cal O}_K$,
the projection
$h\colon W=X\times_YV\to X$ is ${\cal F}$-transversal
on a neighborhood of $h^{-1}(Z)$.

{\rm 1.}
Let
\begin{equation}
\begin{CD}
X@<{p'}<< X'@<{j'}<<W\\
@VfVV @VV{f'}V@VV{f'_V}V\\
Y@<p<<Y'@<j<<V
\end{CD}
\label{eqXYUpi}
\end{equation}
be a cartesian diagram of
regular schemes
of finite type over ${\cal O}_K$.
Assume  that
$p$ is proper
and that $j\colon V=Y'\sm D
\to Y'$ is the open immersion
of the complement of
a divisor $D$ with simple normal crossings.
Then, 
the composition
\begin{equation}
\begin{CD}
{\cal F}\otimes f^*R(pj)_*\Lambda
\to 
{\cal F}\otimes R(p'j')_*\Lambda
@>{\rm (\ref{eqpr0})}>>
R(p'j')_*(p'j')^*{\cal F}
\end{CD}
\label{eqYU}
\end{equation}
where the first morphism is
induced by the base change morphism
is an isomorphism
on a neighborhood of $Z$.

{\rm 2.}
$f$ is universally ${\cal F}$-acyclic
along $Z$.
\end{pr}

For the sake of completeness,
we record the proof in \cite{CC}
with more detail.

\proof{
1. 
Let 
$D_1,\ldots,D_n$
be the irreducible components
of $D$. For a subset $I\subset \{1,\ldots,n\}$,
let $X'_I=X'\times_{Y'}(\bigcap_{i\in I}D_i)$
and let $i'_I\colon X'_I\to X'$
be the closed immersion.
By the assumption,
$p'\colon X'\to X$
and $p'i'_I\colon X'_I\to X$
are ${\cal F}$-transversal
on neighborhoods of
the inverse images of $Z$.

Let ${\cal F}'=p'^*{\cal F}$.
Since the assumption
on $Rh^!\Lambda$
in Proposition \ref{prhF}.1
is satisfied by the absolute
purity \cite[{\sc Th\'eor\`eme 3.1.1}]{purete},
the immersions 
$i'_I\colon X'_I\to X'$
are ${\cal F}'$-transversal
on neighborhoods of
the inverse images of $Z$
by Proposition \ref{prhF}.1.
Hence by Lemma \ref{lmiZ},
the canonical morphism
${\cal F}'
\otimes Rj'_*\Lambda
\to Rj'_*j'^*{\cal F}'$ (\ref{eqpr0}) is an
isomorphism 
on a neighborhood of $p'^{-1}(Z)$.
Since $p'$ is proper,
we obtain an isomorphism
$Rp'_*({\cal F}'
\otimes Rj'_*\Lambda)
\to R(pj')_*(pj')^*{\cal F}$
on a neighborhood of $Z$.

By the projection formula
(\ref{eqprj}),
we have a canonical isomorphism
${\cal F}
\otimes Rp'_*Rj'_*\Lambda
\to Rp'_*({\cal F}'
\otimes Rj'_*\Lambda)$.
The base change morphism
$f^*R(pj)_*\Lambda\to
Rp'_*Rj'_*\Lambda$ is an isomorphism
by the smooth base change
theorem
\cite[Corollaire 1.2]{smbc}.
Hence the morphism (\ref{eqYU})
is an isomorphism on a neighborhood of $Z$.

{\rm 2.}
It suffices to show that
for a smooth morphism
$Y'\to Y$,
the base change
$X'\to Y'$ of $f$
is locally acyclic with respect to the 
pull-back of ${\cal F}$ by Lemma \ref{lmlac}.3.
Similarly as in the proof of 1.,
the assumption is satisfied
for the pull-back $Y'\to Y$.
Hence, 
by replacing $Y$ by $Y'$,
it suffices to show
that
$f$ is locally acyclic with respect to ${\cal F}$.

Let $s\gets t$ be a specialization
of geometric points of $Y$
as in Lemma \ref{lmlac}.1
and let the notation be as loc.~cit.
By \cite[Theorem 4.1,
Theorem 8.2]{dJ},
we may write $t$ as a limit
$\varprojlim_\lambda
U_\lambda$
of the complements $U_\lambda=Y_\lambda\sm 
D_\lambda$,
in regular schemes $Y_\lambda$
endowed with a
proper, surjective
and generically finite
morphism $p_\lambda\colon Y_\lambda
\to Y$ 
of divisors $D_\lambda\subset
Y_\lambda$ with simple normal crossings.
Then, as the limit of
(\ref{eqYU}), the canonical morphism
\begin{equation}
{\cal F}\otimes f_{(s)}^*
Rj_{t*}j^*_t\Lambda
\to 
Rj'_{t*}j^{\prime*}_t{\cal F}
\label{eqijst2}
\end{equation}
is an isomorphism
on the inverse image of $Z$.
Since $Y$ is normal,
the canonical morphism
$\Lambda\to
i_s^*Rj_{t*}j^*_t\Lambda$
is an isomorphism.
Hence the isomorphism
(\ref{eqijst2})
induces an isomorphism
(\ref{eqijst})
on the inverse image of $Z$.
\qed

}

\begin{cor}\label{corlc}
Let $X$
be a regular scheme
of finite type over 
a discrete valuation ring
${\cal O}_K$
and
$Z\subset X$ be a closed subset.
Let
${\cal F}$ be a constructible
complex of $\Lambda$-modules
on $X$.
Assume that every separated morphism
$h\colon W\to X$ of
regular schemes
of finite type over ${\cal O}_K$
is ${\cal F}$-transversal
on a neighborhood of
the inverse image $h^{-1}(Z)$.
Then ${\cal F}$ is locally constant
on a neighborhood of $Z$.
\end{cor}

\proof{
By Proposition \ref{prla1}
applied to $1_X\colon X\to X$,
the identity $1_X\colon X\to X$
is ${\cal F}$-acyclic
along $Z$.
Hence ${\cal F}$ is locally constant
on a neighborhood of $Z$
by Lemma \ref{lmla}.2.
\qed

}

\medskip

We have a partial converse of
Proposition \ref{prla1}
not used in the article.

\begin{pr}[{\cite[Corollary 8.10]{CC}}]\label{prla2}
Let $f\colon X\to Y$
be a smooth morphism
of noetherian schemes
and let
${\cal F}$ be a constructible
complex of $\Lambda$-modules
on $X$.
Let $i\colon Z\to Y$ be an immersion
and let
$$\begin{CD}
X@<h<< W\\
@VfVV@VVV\\
Y@<i<< Z
\end{CD}$$
be a cartesian diagram.
If $f\colon X\to Y$ is ${\cal F}$-acyclic,
then $h\colon W\to X$ is 
${\cal F}$-transversal.
\end{pr}

For the sake of convenience,
we record the proof in \cite{CC}.

\proof{
We may assume that $i\colon Z\to Y$
is a closed immersion.
Let $V=Y\sm Z$
and consider the cartesian diagram
\begin{equation}
\begin{CD}
W@>h>>X@<{j'}<< U\\
@VgVV @VfVV @VV{f_V}V\\
Z@>i>>Y@<j<<V.
\end{CD}
\label{eqZYV}
\end{equation}
By  \cite[Proposition 2.10]{app} applied
to the right square,
we obtain an isomorphism
${\cal F}\otimes f^*Rj_*\Lambda
\to 
Rj'_*j^{\prime*}{\cal F}$.
Since $f$ is smooth,
this induces an isomorphism
${\cal F}\otimes Rj'_*\Lambda
\to 
Rj'_*j^{\prime*}{\cal F}$
by smooth base change theorem
\cite[Corollaire 1.2]{smbc}.
Hence the assertion follows
by Lemma \ref{lmiZ}.2.
\qed

}

\begin{cor}\label{corlabc}
Let $$\begin{CD}
V'@>{g'}>>X'@>{f'}>>Y'\\
@V{h_V}VV@VhVV@VV{h'}V\\
V@>g>> X@>f>>Y
\end{CD}$$
be a cartesian diagram of
morphisms of finite type of schemes
such that 
$f\colon X\to Y$ is smooth
and that the vertical arrows are
separated.
Assume that
$Rh^!\Lambda$ is locally constant
of support $X'$
and that the base change
morphism
$g'^*Rh^!\Lambda
\to Rh_V^!\Lambda$
is an isomorphism.

Let ${\cal G}$
be a constructible complex
of $\Lambda$-modules on $V$
and 
assume that $f$
is $Rg_*{\cal G}$-acyclic
and that $fg$ is
${\cal G}$-acyclic.
Then, 
the base change morphism
\begin{equation}
h^*Rg_*{\cal G}
\to 
Rg'_*h_V^*{\cal G}
\label{eqbcj}
\end{equation}
is an isomorphism.
\end{cor}

\proof{
Since $f$ is $Rg_*{\cal G}$-acyclic
and $fg$ is ${\cal G}$-acyclic,
by Proposition \ref{prla2},
$h$ is $Rg_*{\cal G}$-transversal
and
$h_V$ is ${\cal G}$-transversal.
Hence the assertion follows
from Proposition \ref{prhF}.2.
\qed

}

\section{$C$-transversality}\label{sTX}

In this section,
first we define 
the FW-cotangent bundle
of a regular scheme,
as a vector bundle
on the closed subscheme 
defined by $p=0$.
Then, 
we study properties
of morphisms with respect to
its closed conical subsets
corresponding to the transversality
and the local
acyclicity studied in Section \ref{sF}.

First in Section \ref{ssFW},
we recall basic properties
of the sheaf $F\Omega^1_X$ 
of Frobenius-Witt differentials
from \cite{FW}.
In particular if $X$ is regular,
under a certain finiteness condition,
the sheaf $F\Omega^1_X$
is a locally free 
${\cal O}_{X_{{\mathbf F}_p}}$-module
of rank $\dim X$ on
$X_{{\mathbf F}_p}=
X\times_{{\rm Spec}\, {\mathbf Z}}
{\rm Spec}\, {\mathbf F}_p$.
Under this condition,
we define the FW-cotangent bundle
$FT^*X|_{X_{{\mathbf F}_p}}$ on $X_{{\mathbf F}_p}$
as the vector bundle
associated to the locally free
${\cal O}_{X_{{\mathbf F}_p}}$-module
$F\Omega^1_X$.


We study properties
of morphisms with respect to 
a given closed conical subset in
Sections \ref{ssCtr} and \ref{ssCac}.
In Section \ref{ssCtr},
we study the transversality
for morphisms to $X$.
In Section \ref{ssCac},
we study the acyclicity,
which was also called transversality,
for morphisms from $X$.

\subsection{FW-cotangent bundle}
\label{ssFW}

\begin{df}[{\rm \cite[Definition 1.1]{FW}}]\label{dfFW}
Let $p$ be a prime number.

{\rm 1.}
Define a polynomial
$P\in {\mathbf Z}[X,Y]$
by
\begin{equation}
P=
\sum_{i=1}^{p-1}
\dfrac{(p-1)!}{i!(p-i)!}\cdot
X^iY^{p-i}.
\label{eqP}
\end{equation}

{\rm 2.}
Let $A$ be a ring
and $M$ be an $A$-module.
We say that a mapping
$w\colon A\to M$
is an Frobenius-Witt derivation
or FW-derivation for short
if the following condition is
satisfied:
For any $a,b\in A$, we have
\begin{align}
w(a+b)\, &=
w(a)+
w(b)
-P(a,b)
\cdot w(p),
\label{eqadd}\\
w(ab)\, &=
b^p\cdot w(a)+
a^p\cdot w(b).
\label{eqLb}
\end{align}
\end{df}

Definition \ref{dfFW}.2
is essentially the same
as \cite[Definition 2.1.1]{DKRZ}.
We recall some results from \cite{FW}.

\begin{lm}\label{lmOm}
Let $p$ be a prime number and
$A$ be a ring.

{\rm 1.
(\cite[Lemma 2.1.1]{FW})}
There exists a universal pair
of an $A$-module
$F\Omega^1_A$
and an FW-derivation
$w\colon A
\to F\Omega^1_A$.

{\rm 2.
(\cite[Corollary 2.3.1]{FW})}
If $A$ is a ring over ${\mathbf Z}_{(p)}$,
we have $p\cdot F\Omega^1_A=0$.

{\rm 3.
(\cite[Corollary 2.3.2]{FW})}
If $A$ is a ring over ${\mathbf F}_{p}$,
then there exists a canonical
isomorphism $F\Omega^1_A
\to F^*\Omega^1_A
=\Omega^1_A\otimes_AA$
to the tensor product
with respect to the absolute
Frobenius morphism $A\to A$.
\end{lm}

We call $F\Omega^1_A$
the module of FW-differentials of $A$
and $w(a)\in F\Omega^1_A$
the FW-differential of $a\in A$.
For a morphism $A\to B$ of rings,
we have a canonical $B$-linear morphism
$F\Omega^1_A\otimes_AB
\to
F\Omega^1_B$.

We may sheafify the construction
and define $F\Omega^1$
as a quasi-coherent ${\cal O}_X$-module
for a scheme $X$. 
We call $F\Omega^1_X$
the sheaf of FW-differentials on $X$.
If $X$ is a scheme over
${\mathbf Z}_{(p)}$,
the ${\cal O}_X$-module
$F\Omega^1_X$ is
an ${\cal O}_{X_{{\mathbf F}_p}}$-module
where
$X_{{\mathbf F}_p}
=X\times_{{\rm Spec}\, {\mathbf Z}_{(p)}}
{\rm Spec}\, {\mathbf F}_p$.
Further if $X$ is noetherian
and if $X_{{\mathbf F}_p}$
is of finite type over a field
of finite $p$-basis,
then $F\Omega^1_X$ is
a coherent 
${\cal O}_{X_{{\mathbf F}_p}}$-module
by \cite[Lemma 4.1.2]{FW}.
If $X$ is a scheme over
${\mathbf F}_p$,
we have a canonical isomorphism
\begin{equation}
F\Omega^1_X
\to F^*\Omega^1_X
\label{eqFFX}
\end{equation}
to the pull-back by
the absolute Frobenius morphism
$F\colon X\to X$,
sending $w(a)$ to $da$.

For a morphism 
$f\colon X\to Y$ of schemes,
we have a canonical morphism
\begin{equation}
f^*F\Omega^1_Y\to
F\Omega^1_X
\label{eqFXY}
\end{equation}

\begin{pr}[{\rm \cite[Proposition 2.4]{FW}}]\label{prdx}
Let $X$ be a scheme
and $x\in X$
be a point such that
the residue field $k(x)={\cal O}_{X,x}/
{\mathfrak m}_{X,x}$
is of characteristic $p$.
For a $k(x)$-vector space $M$,
let $F^*M$
denote the tensor product
$M\otimes_{k(x)}k(x)$
with respect to the Frobenius
$F\colon k(x)\to k(x)$.
Then,  we have an exact
sequence 
\begin{equation}
\begin{CD}
0@>>>
F^*({\mathfrak m}_{X,x}/
{\mathfrak m}_{X,x}^2)
@>{w}>>
F\Omega^1_{X,x}
\otimes_{{\cal O}_{X,x}} k(x)
@>{\rm (\ref{eqFFX})}>>
F^*\Omega^1_{k(x)}
@>>>0
\end{CD}
\label{eqdx}
\end{equation}
of $k(x)$-vector spaces.
\end{pr}

\begin{pr}[{\rm \cite[Proposition 2.8]{FW}}]\label{prsm}
Let $f\colon X\to Y$
be a morphism of finite type of
regular noetherian schemes
over ${\mathbf Z}_{(p)}$.
Then the following conditions
are equivalent:

{\rm (1)}
$f\colon X\to Y$ is smooth
on a neighborhood of
$X_{{\mathbf F}_p}$.

{\rm (2)}
The sequence 
\begin{equation}
\begin{CD}
0@>>>
f^*F\Omega^1_Y
@>{\rm (\ref{eqFXY})}>>
F\Omega^1_X
@>{\rm (\ref{eqFFX})}>>F^*\Omega^1_{
X_{{\mathbf F}_p}/
Y_{{\mathbf F}_p}}
@>>>
0
\end{CD}
\end{equation}
of ${\cal O}_{X_{{\mathbf F}_p}}$-modules
is a locally split exact sequence.
\end{pr}

\begin{thm}[{\rm \cite[Theorem 3.1]{FW}}]
\label{thmreg}
Let $X$ be a noetherian scheme
over ${\mathbf Z}_{(p)}$
and $X_{{\mathbf F}_p}=
X\times_{{\rm Spec}\, 
{\mathbf Z}_{(p)}}
{\rm Spec}\, {\mathbf F}_p$
be the closed subscheme.
Assume that the reduced
part $X_{{\mathbf F}_p,{\rm red}}$ is
a scheme of finite type over
a field $k$ with finite $p$-basis.
%
%
%
If $X$ is regular and 
is equi-dimensional
of dimension $n$ and if
$[k:k^p]=p^r$, then 
the ${\cal O}_{X_{{\mathbf F}_p}}$-module
$F\Omega^1_X$
is locally free of rank $n+r$.
\end{thm}

\begin{cor}[{\rm \cite[Corollary 2.6,
Corollary 3.2]{FW}}]
\label{corXZ}
Let $X$ be a regular noetherian scheme
over ${\mathbf Z}_{(p)}$
such that the reduced
part $X_{{\mathbf F}_p,{\rm red}}$ 
of
$X_{{\mathbf F}_p}=X\times_{{\rm Spec}\, 
{\mathbf Z}_{(p)}}
{\rm Spec}\, {\mathbf F}_p$
is
a scheme of finite type over
a field $k$ of finite $p$-basis.
Let $Z\subset X$ be a closed subscheme.

We consider the following conditions:

{\rm (1)}
$Z$ is regular on a neighborhood of
$Z_{{\mathbf F}_p}=Z\times_{{\rm Spec}\, 
{\mathbf Z}_{(p)}}
{\rm Spec}\, {\mathbf F}_p$.

{\rm (1$'$)}
At every point $x\in Z
_{{\mathbf F}_p}$,
the local ring
${\cal O}_{Z,x}$ is regular.

{\rm (2)}
The sequence 
\begin{equation}
\begin{CD}
0@>>>
F^*(N_{Z/X}\otimes_{{\cal O}_Z}
{\cal O}_{Z_{{\mathbf F}_p}})
@>w>>
F\Omega^1_X\otimes_{{\cal O}_X}
{\cal O}_{Z_{{\mathbf F}_p}}
\longrightarrow
F\Omega^1_Z
@>>>
0
\end{CD}
\label{eqXZ}
\end{equation}
of 
${\cal O}_{Z_{{\mathbf F}_p}}$-modules
is a locally splitting exact sequence.

Then, we have
{\rm (1)}$\Rightarrow${\rm (2)}$\Rightarrow${\rm (1$'$)}.
Consequently 
if the subset
${\rm Reg}(Z)
\subset Z$ consisting
of regular points is an open subset,
the 3 conditions are equivalent.
\end{cor}

\proof{
The implications
{\rm (1)}$\Rightarrow${\rm (2)} and
{\rm (2)}$\Rightarrow${\rm (1$'$)}
are proved in 
\cite[Corollary 3.2]{FW} and in
\cite[Corollary 2.6.1]{FW}
respectively.
Since (1$'$) means
$Z_{{\mathbf F}_p}
\subset {\rm Reg}(Z)$,
the last assertion follows.
\qed

}


\begin{df}\label{dfFTX}
Let $k$ be a perfect field 
of characteristic $p>0$
and let $X$ be
a regular noetherian scheme
satisfying the following condition:

{\rm (F)}
$X_{{\mathbf F}_p}
=X\times_{{\rm Spec}\, {\mathbf Z}}
{\rm Spec}\, {\mathbf F}_p$
is a scheme of
finite type over $k$.

\noindent
Then, we define
the FW-cotangent bundle
$FT^*X|_{X_{{\mathbf F}_p}}$
of $X$ to be the vector bundle
on $X_{{\mathbf F}_p}$
associated with the locally free
${\cal O}_{X_{{\mathbf F}_p}}$-module
$F\Omega^1_X$
of rank $\dim X$.
\end{df}

Let $x\in X_{{\mathbf F}_p}$
be a closed point 
and let
$T^*_xX$
denote the cotangent space
at $x$ defined as a scheme
${\rm Spec}\,
S_{k(x)}({\mathfrak m}_x/
{\mathfrak m}_x^2)^\vee$
associated to the
$k(x)$-vector space
${\mathfrak m}_x/
{\mathfrak m}_x^2$.
Since $k(x)$ is perfect,
the exact sequence
(\ref{eqdx}) defines a canonical
isomorphism
\begin{equation}
F^*T^*_xX
\to 
FT^*X|_x
\label{eqTx}
\end{equation}
to the fiber of 
the FW-cotangent bundle
at $x$ 
from the pull-back by Frobenius
$F\colon x\to x$ of 
$T^*_xX$.
If $X=X_{{\mathbf F}_p}$,
then the FW-cotangent bundle
$FT^*X|_{X_{{\mathbf F}_p}}$
is the pull-back of
the cotangent bundle
$T^*X$ 
by the Frobenius morphism
$F\colon X\to X$
by {\rm (\ref{eqFFX})}.

Let $X\to Y$ be a morphism
of finite type
of regular noetherian schemes
satisfying the condition (F)
in Definition \ref{dfFTX}.
Then, the morphism (\ref{eqFXY})
defines morphisms
\begin{equation}
\begin{CD}
FT^*X|_{X_{{\mathbf F}_p}}
@<{f^*}<<
FT^*Y|_{Y_{{\mathbf F}_p}}
\times_{Y_{{\mathbf F}_p}}
X_{{\mathbf F}_p}
@>>>
FT^*Y|_{Y_{{\mathbf F}_p}}
\end{CD}
\label{eqdXY}
\end{equation}
 of
schemes.

Assume that $X\to Y$ is
smooth and let 
$F^*T^*X/Y|_{
X_{{\mathbf F}_p}}$ denote the
pull-back by the Frobenius
$F\colon X_{{\mathbf F}_p}
\to X_{{\mathbf F}_p}$
of the restriction to
$X_{{\mathbf F}_p}$ of 
the vector
bundle defined $T^*X/Y$ by
the locally free ${\cal O}_X$-module
$\Omega^1_{X/Y}$.
Then, by Proposition \ref{prsm},
we have an exact sequence
\begin{equation}
0\to FT^*Y|_{Y_{{\mathbf F}_p}}
\times_{Y_{{\mathbf F}_p}}
X_{{\mathbf F}_p}
\to FT^*X|_{X_{{\mathbf F}_p}}
\to F^*T^*X/Y|_{
X_{{\mathbf F}_p}}\to 0
\label{eqTEf}
\end{equation}
of vector bundles on $X_{{\mathbf F}_p}$.

Similarly,
let $Z\to X$ be a closed immersion
of regular noetherian schemes
satisfying the condition (F).
Let ${\cal I}_Z\subset {\cal O}_X$
be the ideal sheaf
and let $T^*_ZX$ be
the conormal bundle
defined by
the locally free ${\cal O}_Z$-module
${\cal I}_Z/{\cal I}_Z^2$.
Let 
$F^*T^*_ZX|_{
Z_{{\mathbf F}_p}}$ denote the
pull-back by the Frobenius
$F\colon Z_{{\mathbf F}_p}
\to Z_{{\mathbf F}_p}$
of the restriction to
$Z_{{\mathbf F}_p}$.
Then, by Corollary \ref{corXZ},
we have an exact sequence
\begin{equation}
0\to F^*T^*_ZX|_{
Z_{{\mathbf F}_p}}\to 
FT^*X|_{Z_{{\mathbf F}_p}}
\to 
FT^*Z|_{Z_{{\mathbf F}_p}}
\to 0
\label{eqTEi}
\end{equation}
of vector bundles on $Z_{{\mathbf F}_p}$.

\subsection{$C$-transversality}\label{ssCtr}

In the rest of this section,
we fix a perfect field $k$
of characteristic $p>0$.

We fix some terminology
on closed conical subsets of
a vector bundle of a scheme.
Let $V$ be a vector bundle
over a scheme $Y$.
We say that a closed subset
of $V$ is conical if it is
stable under the action of
${\mathbf G}_{m,Y}$.
For a closed conical subset
$C\subset V$,
the intersection
$B=C\cap Y$ with the
$0$-section $Y\subset V$ regarded
as a closed subset of $Y$
is called the base of $C$.
The base $B$ equals the
image of $C$ by
the projection $V\to Y$.

We say that a separated
morphism $f\colon X\to Y$
of finite type of schemes
is proper on a closed subset $Z\subset X$
if for every base change
$f'\colon X'\to Y'$ of $f$
its restriction to
the inverse image $Z'
\subset X'$ is a closed mapping.
For a morphism
$V\to V'$ of vector bundles
on a scheme $Y$
and a closed conical subset
$C$ of $V$,
the morphism $V\to V'$
is proper on $C$ if and only
if the intersection
$C\cap {\rm Ker}(V\to V')$
is a subset of the $0$-section of $V$
by \cite[Lemma 1.2(ii)]{Be}.

\begin{df}\label{dfhC}
Let $X$ be a 
regular noetherian scheme
satisfying the condition {\rm (F)}
in Definition {\rm \ref{dfFTX}}
and let $C\subset FT^*X|_{X_{{\mathbf F}_p}}
$ be a closed
conical subset
of the FW-cotangent bundle.
Let $h\colon W\to X$
be a morphism of finite type
of regular schemes.

{\rm 1.}
{\rm (\cite[1.2]{Be}, \cite[Definition 3.3]{CC})}
We say that 
$h\colon W\to X$ is $C$-transversal
if 
the intersection of
$h^*C=C\times_XW
\subset 
FT^*X|_{X_{{\mathbf F}_p}}
\times_{X_{{\mathbf F}_p}}
W_{{\mathbf F}_p}$
with the kernel 
${\rm Ker}(
FT^*X|_{X_{{\mathbf F}_p}}
\times_{X_{{\mathbf F}_p}}
W_{{\mathbf F}_p}
\to
FT^*W|_{W_{{\mathbf F}_p}})$
is a subset of the $0$-section.

{\rm 2.}
Assume that $h$ is
$C$-transversal.
Then we define
a closed conical subset 
$h^\circ C
\subset FT^*W|_{W_{{\mathbf F}_p}}$
to be the image
of $h^*C$ by
$
FT^*X|_{X_{{\mathbf F}_p}}
\times_{X_{{\mathbf F}_p}}
W_{{\mathbf F}_p}
\to
FT^*W|_{W_{{\mathbf F}_p}}$.
\end{df}

Example.
Let $Z\subset X$ be a regular closed subscheme.
Then a closed
conical subset 
$C=F^*T^*_ZX|_{Z_{{\mathbf F}_p}}
\subset 
FT^*X|_{X_{{\mathbf F}_p}}$
is defined by (\ref{eqTEi}).
In particular,
for $Z=X$,
the $0$-section
$F^*T^*_XX|_{X_{{\mathbf F}_p}}
=X_{{\mathbf F}_p}$ is
a closed conical subset of
$FT^*X|_{X_{{\mathbf F}_p}}$.

\begin{lm}\label{lmTXC}
Let $X$ be a 
regular noetherian scheme
satisfying the condition {\rm (F)}
in Definition {\rm \ref{dfFTX}}
and let $C\subset FT^*X|_{X_{{\mathbf F}_p}}
$ be a closed
conical subset.
Let $h\colon W\to X$
be a morphism of finite type
of regular schemes.

{\rm 1.}
Let $C=FT^*X|_Z$
be the restriction to 
a closed subset
$Z\subset X_{{\mathbf F}_p}$
of the closed fiber.
If $h$ is $C$-transversal,
then $h$ is smooth
on a neighborhood of
the inverse image $h^{-1}(Z)$.

{\rm 2.}
If $C$ is the $0$-section
of $FT^*X|_{X_{{\mathbf F}_p}}$,
then $h$ is $C$-transversal.

{\rm 3.}
If $h$ is smooth,
for any closed conical subset
$C$ of $FT^*X|_{X_{{\mathbf F}_p}}$,
the morphism
$h$ is $C$-transversal.
\end{lm}

\proof{
1.
The condition that the
intersection of
$h^*C=FT^*X|_Z
\times_{X_{{\mathbf F}_p}}
W_{{\mathbf F}_p}
=FT^*X
\times_{X_{{\mathbf F}_p}}
h^{-1}(Z)$
with the kernel
${\rm Ker}(FT^*X|_{X_{{\mathbf F}_p}}
\times_{X_{{\mathbf F}_p}}
W_{{\mathbf F}_p}
\to FT^*W|_{W_{{\mathbf F}_p}})$
is a subset of the $0$-section
means that
$F\Omega^1_X
\otimes_{{\cal O}_X}{\cal O}_W
\to F\Omega^1_W$
is a locally splitting injection
on a neighborhood of $h^{-1}(Z)$.
By Proposition \ref{prsm},
this means that
$W\to X$
is smooth on a neighborhood of
the inverse image $h^{-1}(Z)$.

2.
If $C$ is the $0$-section,
its intersection with the
kernel 
${\rm Ker}(FT^*X|_{X_{{\mathbf F}_p}}
\times_{X_{{\mathbf F}_p}}
W_{{\mathbf F}_p}
\to FT^*W|_{W_{{\mathbf F}_p}})$
is also the $0$-section.

3.
If $h$ is smooth,
the morphism $FT^*X|_{X_{{\mathbf F}_p}}
\times_{X_{{\mathbf F}_p}}
W_{{\mathbf F}_p}
\to FT^*W|_{W_{{\mathbf F}_p}}$ is an injection
by Proposition \ref{prsm}.
Hence 
for any subset
$C\subset FT^*X|_{X_{{\mathbf F}_p}}$,
its intersection with the
kernel 
${\rm Ker}(FT^*X|_{X_{{\mathbf F}_p}}
\times_{X_{{\mathbf F}_p}}
W_{{\mathbf F}_p}
\to FT^*W|_{W_{{\mathbf F}_p}})$
is a subset of the $0$-section.
\qed

}

\begin{lm}\label{lmhC}
Let $h\colon W\to X$
be a morphism of
finite type of 
regular noetherian schemes
satisfying the condition {\rm (F)}
and let $C$ be a closed
conical subset of $FT^*X|_{X_{{\mathbf F}_p}}$.
Assume that $h$ is
$C$-transversal.
Then, for a morphism 
$g\colon V\to W$ of finite type of
regular noetherian schemes
the following conditions
are equivalent:

{\rm (1)}
The morphism
$g$ is $h^\circ C$-transversal.

{\rm (2)}
The composition
$hg$ is $C$-transversal.

\noindent
If these equivalent conditions
are satisfied,
we have $(hg)^\circ C=g^\circ h^\circ C$.
\end{lm}

\proof{
The condition (1)
means that
the intersection
$h^*C\cap {\rm Ker}(
FT^*X|_{X_{{\mathbf F}_p}}
\times_{X_{{\mathbf F}_p}}
W_{{\mathbf F}_p}
\to FT^*W|_{W_{{\mathbf F}_p}})$
is a subset of the $0$-section
and further that
for $h^\circ C
\subset FT^*W|_{W_{{\mathbf F}_p}}$,
the intersection
$g^*h^\circ C\cap {\rm Ker}(
FT^*W|_{W_{{\mathbf F}_p}}
\times_{W_{{\mathbf F}_p}}
V_{{\mathbf F}_p}
\to FT^*V|_{V_{{\mathbf F}_p}})$
is a subset of the $0$-section.
This means that
$(hg)^*C\cap {\rm Ker}(
FT^*X|_{X_{{\mathbf F}_p}}
\times_{X_{{\mathbf F}_p}}
V_{{\mathbf F}_p}
\to FT^*V|_{V_{{\mathbf F}_p}})$
is a subset of the $0$-section,
namely the condition (2).

The image of 
$(hg)^*C$ by
$FT^*X|_{X_{{\mathbf F}_p}}
\times_{X_{{\mathbf F}_p}}
V_{{\mathbf F}_p}
\to FT^*V|_{V_{{\mathbf F}_p}}$
equals 
the image of $g^*h^\circ C$
by 
$FT^*W|_{W_{{\mathbf F}_p}}
\times_{W_{{\mathbf F}_p}}
V_{{\mathbf F}_p}
\to FT^*V|_{V_{{\mathbf F}_p}}$.
\qed

}

\medskip
The terminology transversality
is related to the transversality
of morphisms of regular
schemes defined as follows.

\begin{df}\label{dftrans}
Let $f\colon X\to Y$
and $g\colon V\to Y$
be morphisms of finite type of
regular schemes
and set $W
=X\times_YV$.

{\rm 1.}
Let $w\in W$
and $x\in X,y\in Y,v\in V$
be the images.
We say that 
$f$ and $g$ are transversal
at $w$, if ${\cal O}_{W,w}$
is regular and if
$Tor_q^{{\cal O}_{Y,y}}
({\cal O}_{X,x},{\cal O}_{V,v})=0$
for $q>0$.

{\rm 2.}
Let $W_1\subset W$
be an open subscheme.
We say that 
$f$ and $g$ are transversal on
$W_1$ if 
$f$ and $g$ are transversal 
at every point of $W_1$.
\end{df}

Example.
Let $Z\subset X$ be a regular closed subscheme
and 
$C=F^*T^*_ZX|_{Z_{{\mathbf F}_p}}
\subset 
FT^*X|_{Z_{{\mathbf F}_p}}$
be the closed
conical subset defined by the conormal bundle.
Then, as we will see
in Corollary \ref{corfC},
a morphism
$h\colon W\to X$ of
finite type 
of regular quasi-excellent
noetherian schemes
is $C$-transversal
if and only if
$h\colon W\to X$
is transversal to $Z\subset X$
on a neighborhood of
the closed fiber $W_{{\mathbf F}_p}$.

In particular, 
if $X$ is smooth over 
a discrete valuation 
ring ${\cal O}_K$ of mixed
characteristic with residue field $k$
and  if $C=F^*T^*_{X_k}X|_{X_k}$
for the closed fiber $Z=X_k$,
then the condition that
$h\colon W\to X$ is
$C$-transversal
means that
$W$ is smooth over ${\cal O}_K$ 
on a neighborhood of
the closed fiber $W_k$.


\begin{lm}\label{lmtrreg}
Let $f\colon X\to Y$
and $g\colon V\to Y$
be morphisms of finite type of
regular schemes
and set $W
=X\times_YV$.
Let $w\in W$
and $x\in X,y\in Y,v\in V$
be the images.

{\rm 1.}
Suppose that $g\colon V\to Y$
is an immersion.
Then, the following conditions
are equivalent:

{\rm (1)}
$f$ and $g$ are transversal
at $w$.

{\rm (2)}
The morphism
$T^*_yY\times_yx\to T^*_xX$
on the cotangent space
induces an injection
on the subspace
$T^*_VY\times_Vy
\subset
T^*_yY\times_yx$.

{\rm 2.}
Suppose that the
subset ${\rm Reg}(W)
\subset W$
consisting of regular points
is an open subset.
If
$f$ and $g$ are transversal
at $w\in W$,
then
$f$ and $g$ are transversal
on a neighborhood of $w$.
\end{lm}

The condition that
${\rm Reg}(W)
\subset W$
is an open subset is satisfied
if $W$ is of finite type
over a Dedekind domain
such that the fraction field
is of characteristic $0$
or a semi-local ring of dimension
at most $1$
by \cite[Corollaire (6.12.6)]{EGA4}.

\proof{
1.
Let $a_1,\ldots,a_r\in {\cal O}_{Y,y}$
be a minimal system of generators
of ${\rm Ker}({\cal O}_{Y,y}
\to {\cal O}_{V,y})$.
Then, the both conditions are
equivalent to the condition
that $a_1,\ldots,a_r\in {\cal O}_{X,x}$
is a part of a regular system of
parameters.

2.
Since the ${\cal O}_W$-modules
${\cal T}or_q^{{\cal O}_Y}
({\cal O}_X,{\cal O}_V)=0$
are coherent
and $w$ is an element
of the open subset ${\rm Reg}(W)$,
the assertion follows.
\qed

}

\medskip

Let $f\colon X\to Y$
be a morphism of finite type of
regular noetherian schemes
such that $Y_{{\mathbf F}_p}$
is of finite type
over $k$
and consider the
morphisms
(\ref{eqdXY}).
Let $C$ be a closed conical subset
of $FT^*X|_{X_{{\mathbf F}_p}}$
such that 
$f\colon X\to Y$ is proper
on the base $B(C)$.
Then we define a closed conical subset
$f_\circ C$ of $FT^*Y|_{Y_{{\mathbf F}_p}}$
to be
the image by
$FT^*X|_{X_{{\mathbf F}_p}}
\times_{X_{{\mathbf F}_p}}
Y_{{\mathbf F}_p}
\to
FT^*Y|_{Y_{{\mathbf F}_p}}$
of the inverse image of
$C$ by
$FT^*X|_{X_{{\mathbf F}_p}}
\times_{X_{{\mathbf F}_p}}
Y_{{\mathbf F}_p}
\to
FT^*X|_{X_{{\mathbf F}_p}}$.

For a closed immersion
$i\colon Z\to X$
of regular noetherian schemes
such that $X_{{\mathbf F}_p}$
is of finite type
over $k$,
the closed conical subset
$F^*T^*_ZX|_{X_{{\mathbf F}_p}}$
defined by the conormal bundle
equals $i_\circ C$
for the $0$-section
$C=FT^*_ZZ|_{Z_{{\mathbf F}_p}}$
of $FT^*Z|_{Z_{{\mathbf F}_p}}$.

\begin{pr}\label{prfC}
Let $X,Y$ and $V$ be regular
noetherian schemes 
satisfying the condition {\rm (F)} and
\begin{equation}
\begin{CD}
X@<h<< W\\
@VfVV@VV{f'}V\\
Y@<g<< V
\end{CD}
\label{eqprfC}
\end{equation}
be a cartesian
diagram of morphisms of
finite type.
Assume that the
subset ${\rm Reg}(W)
\subset W$
consisting of regular points
is an open subset.
Let $C$ be a closed
conical subset of $FT^*X|_{X_{{\mathbf F}_p}}$
such that $f$ is proper on the
base $B(C)$.
Then,
the following conditions
are equivalent:

{\rm (1)}
The morphism
$g$ is $f_\circ C$-transversal.

{\rm (2)}
There exists a
regular neighborhood $W_1
\subset W$  
of the inverse image of
the base $B(C)$
such that 
$f$ and $g$ are transversal on
$W_1$ and that the restriction
$h_1\colon W_1\to X$ is $C$-transversal.

\noindent
If these equivalent conditions are satisfied,
we have $g^\circ f_\circ C=
f'_{1\circ} h_1^\circ C$
for the restriction
$f'_1\colon W_1\to V$ of $f'$.
\end{pr}

\proof{
(1)$\Rightarrow$(2):
Let $x\in B(C)\subset X_{{\mathbf F}_p}$
be a closed point and $y=f(x)\in Y_{{\mathbf F}_p}$.
Since the assertion is \'etale local,
we may assume that
the morphism $k(y)\to k(x)$
of residue fields is an isomorphism.
There exist
an open neighborhood 
$U\subset X$ of $x\in X$
and a cartesian diagram
$$\begin{CD}
W@<{\supset}<<W\times_XU
@>>> Q@>>> V\\
@VVV@VVV@VVV@VVV\\
X@<{\supset}<<U@>>> P@>>>Y
\end{CD}$$
such that $P\to Y$
is smooth and 
$U\to P$ is a closed immersion.
Let $w\in W$ be a closed point
above $x$ and $v=f'(w)\in V$.
We may also assume that
the morphisms 
$k(y)\to k(v)$
and hence $k(x)\to k(w)$
are isomorphisms.
We consider the cartesian
diagram
$$\begin{CD}
@.
T^*_wQ
@<<<
T^*_vV
\\
@.@AAA@AAA\\
T^*_xX
@<<<
T^*_xP
@<<<
T^*_yY
\end{CD}$$
of cotangent spaces
and identify their Frobenius
pull-backs with the fibers
of FW-cotangent bundles
by the isomorphism (\ref{eqTx}).

Let $\widetilde C_x
\subset F^*T^*_xP$
and $A_x \subset F^*T^*_yY$
be the inverse images of
$C_x\subset F^*T^*_xX$.
Then, by the condition (1),
the intersection 
$A_x\cap {\rm Ker}(F^*T^*_yY \to 
F^*T^*_vV)$
is a subset of the $0$-section.
Since
$T^*_yY \to T^*_xP$
induces an isomorphism
${\rm Ker}(T^*_yY\to T^*_vV)
\to 
{\rm Ker}(T^*_xP\to T^*_wQ)$,
the intersection 
$\widetilde C_x\cap
{\rm Ker}(F^*T^*_xP\to F^*T^*_wQ)$
is a subset of the $0$-section.

By the exact sequence
$0\to T^*_XP|_x
\to T^*_xP\to T^*_xX\to 0$
and $x\in B(C)$,
we have
$F^*T^*_XP|_x\subset \widetilde C_x$.
Hence
$T^*_xP\to T^*_wQ$
induces an injection on
$T^*_XP|_x$.
Namely, 
the morphism
$Q\to P$ and the immersion
$U\to P$ are transversal
on a neighborhood of $w$
by Lemma \ref{lmtrreg}.

Hence
the horizontal arrows
of the commutative diagram
\begin{equation}
\begin{CD}
T^*_wW
@<<<
T^*_vV
\\
@AAA@AAA\\
T^*_xX
@<<<
T^*_yY
\end{CD}
\label{eqUVW}
\end{equation}
induce isomorphisms
on the kernels and cokernels
of the vertical arrows.
Since the intersection of
the inverse image $A_x$ with 
${\rm Ker}(F^*T^*_yY
\to F^*T^*_wV)$
is a subset of the $0$-section,
the intersection of
$C_x$ with 
${\rm Ker}(F^*T^*_xX
\to F^*T^*_wW)$
is also a subset of the $0$-section.
Namely,
$h$ is $C$-transversal
on a neighborhood of $w$.
Thus $h$ is $C$-transversal on
a neighborhood of 
the inverse image of $B(C)$.

Further an elementary
diagram chasing shows
that the inverse image of
$h^\circ C|_w$
by $F^*T^*_wW\gets F^*T^*_vV$
equals the image of
$A_x$ by
$F^*T^*_yY\to F^*T^*_vV$.
Hence we have
$g^\circ f_\circ C=
f'_{1\circ} h_1^\circ C$.

(2)$\Rightarrow$(1):
Let $w\in B(h_1^\circ C)$
be a closed point
and let $v\in V, x\in X$
and $y\in Y$ be the image.
Then, the commutative diagram
(\ref{eqUVW})
induces an isomorphism
${\rm Ker}(T^*_yY\to T^*_vV)
\to
{\rm Ker}(T^*_xX\to T^*_wW)$
on the kernels.
In the same notation,
since the intersection of
$C_x$ with 
${\rm Ker}(F^*T^*_xX\to F^*T^*_wW)$
is a subset of the $0$-section,
the intersection of
$A_x$ with 
${\rm Ker}(F^*T^*_yY\to F^*T^*_vV)$
is also a subset of the $0$-section.
\qed

}

\begin{cor}\label{corfC}
Let $X,Y$ and $V$
be regular noetherian schemes 
satisfying the condition {\rm (F)}
and let
{\rm (\ref{eqprfC})}
be a cartesian diagram
of morphisms of finite type.
Assume that the
subset ${\rm Reg}(W)
\subset W$
consisting of regular points
is an open subset
and that
$f\colon X\to Y$ is an immersion.
Then,
the following conditions
are equivalent:

{\rm (1)}
The morphism
$g$ is $F^*T^*_XY|_{Y_{{\mathbf F}_p}}$-transversal.

{\rm (2)}
The morphism $g\colon V\to Y$
is transversal with
the immersion
$f\colon X\to Y$ on
a neighborhood of $W_{{\mathbf F}_p}
=V\times_XX_{{\mathbf F}_p}$.
\end{cor}

\proof{
It suffices to apply
Proposition \ref{prfC} 
together with Lemma \ref{lmTXC}.2
to
the $0$-section $C$
of $FT^*X|_{X_{{\mathbf F}_p}}$.
\qed

}

\begin{df}\label{dfCet}
Let $f\colon U\to X$ be an \'etale morphism 
of regular noetherian schemes 
satisfying the condition {\rm (F)}
and let $C'$ be a closed
conical subset of $FT^*U$.
We identify $FT^*U$ with
the pull-back
$FT^*X\times_{X_{{\mathbf F}_p}}
U_{{\mathbf F}_p}$ by the
canonical isomorphism
induced by
$F\Omega^1_X\otimes_{{\cal O}_X}
{\cal O}_U\to
F\Omega^1_U$ 
and let
${\rm pr}_1\colon 
FT^*X\times_{X_{{\mathbf F}_p}}
U_{{\mathbf F}_p} \to 
FT^*X$ be the projection.
Then, we define
a closed conical subset $f_*C'$ of 
$FT^*X$ to be the union of
the closure $\overline{{\rm pr}_1(C')}$
and the restriction
$FT^*X|_{X_{{\mathbf F}_p}\sm 
f(U_{{\mathbf F}_p})}$ 
to the image of the complement.
\end{df}

\begin{lm}\label{lmCet}
Let $$
\begin{CD}
V@>>>W\\
@VgVV@VVhV\\
U@>f>>X
\end{CD}
$$ be a cartesian diagram
of regular noetherian schemes 
satisfying the condition {\rm (F)}
such that $f$ is
an \'etale morphism of finite type.
Let $C'$ be a closed
conical subset of $FT^*U$
and set $C=f_*C'\subset FT^*X$
as in Definition {\rm \ref{dfCet}}.

If $h$ is $C$-transversal,
then $h$ is smooth on a neighborhood
of $h^{-1}(X_{{\mathbf F}_p}\sm 
f(U_{{\mathbf F}_p}))$
and $g$ is $C'$-transversal.
\end{lm}

\proof{
Assume that
$h$ is $C$-transversal.
Since 
$C\supset 
FT^*X|_{X_{{\mathbf F}_p}\sm 
h(U_{{\mathbf F}_p})}$,
the morphism 
$h$ is smooth on a neighborhood
of $h^{-1}(X_{{\mathbf F}_p}\sm 
f(W_{{\mathbf F}_p}))$
by Lemma \ref{lmTXC}.1.
Since 
$f^\circ C\supset C'$,
the morphism 
$g$ is $C'$-transversal
by Lemma \ref{lmhC}.
\qed

}

\subsection{$C$-acyclicity}\label{ssCac}

We keep fixing a perfect field
$k$ of characteristic $p>0$.

\begin{df}\label{dffC}
Let $f\colon X\to Y$
be a morphism of finite type of 
regular noetherian schemes
satisfying the condition {\rm (F)}
in Definition {\rm \ref{dfFTX}}
and 
let $C$
be a closed conical subset
of the FW-cotangent bundle
$FT^*X|_{X_{{\mathbf F}_p}}$.
We say that $f$
is $C$-acyclic if the inverse image of
$C$ by the morphism
$FT^*Y|_{Y_{{\mathbf F}_p}}
\times_{Y_{{\mathbf F}_p}}
X_{{\mathbf F}_p}
\to FT^*X|_{{\mathbf F}_p}$
is a subset of the $0$-section.
\end{df}


The corresponding notion is
called $C$-transversality
in \cite[1.2]{Be} and
\cite[Definition 3.5]{CC}.
Here to avoid confusion with
the $C$-transversality 
for morphisms to $X$ in
Definition \ref{dfhC}.1,
\cite[1.2]{Be} and
\cite[Definition 3.3]{CC},
we introduce another terminology.
We will show in Lemma \ref{lmhf}.2
that
for a morphism
$f\colon X\to Y$
of regular schemes and
a closed immersion
$i\colon Z\to X$
of regular schemes,
the morphism
$f$ is $F^*T^*_ZX|_{X_{{\mathbf F}_p}}$-acyclic
if and only if
the composition
$fi$ is smooth on a neighborhood
of $Z_{{\mathbf F}_p}$.

\begin{lm}\label{lmftr}
Let $f\colon X\to Y$
be a morphism of finite type of 
regular noetherian schemes
satisfying the condition {\rm (F)}
and 
let $C$
be a closed conical subset
of $FT^*X|_{X_{{\mathbf F}_p}}$.

{\rm 1.}
The following conditions
are equivalent:

{\rm (1)}
$f$ is $C$-acyclic.

{\rm (2)}
$f$ is smooth
on a neighborhood of
the base $B(C)\subset 
X_{{\mathbf F}_p}$
and 
the intersection of
$C\subset FT^*X|_{X_{{\mathbf F}_p}}$ 
with the image of the morphism
$FT^*Y|_{Y_{{\mathbf F}_p}}
\times_{Y_{{\mathbf F}_p}}X_{{\mathbf F}_p}
\to FT^*X|_{X_{{\mathbf F}_p}}$
is a subset of the $0$-section.

{\rm 2.}
If $C$ is the $0$-section
$F^*T^*_XX|_{X_{{\mathbf F}_p}}$,
the following conditions
are equivalent:

{\rm (1)}
$f$ is $C$-acyclic.

{\rm (2)}
$f$ is smooth
on the neighborhood of
$X_{{\mathbf F}_p}$.
\end{lm}

\proof{
1.
The condition (1)
is equivalent to the conjunction
of the following (1$'$) and (1$''$):

(1$'$)
The inverse image of
the $0$-section by
$FT^*Y|_{Y_{{\mathbf F}_p}}
\times_{Y_{{\mathbf F}_p}}X_{{\mathbf F}_p}
\to FT^*X|_{X_{{\mathbf F}_p}}$
on the base $B(C)\subset X_{{\mathbf F}_p}$
is a subset of the $0$-sections.

(1$''$)
The intersection of
$C\subset FT^*X|_{X_{{\mathbf F}_p}}$ 
with the image of the morphism
$FT^*Y|_{Y_{{\mathbf F}_p}}
\times_{Y_{{\mathbf F}_p}}X_{{\mathbf F}_p}
\to FT^*X|_{X_{{\mathbf F}_p}}$
is a subset of the $0$-sections.

\noindent
The condition (1$'$)
means that
the morphism
$f^*F\Omega^1_Y
\to
F\Omega^1_X$
is a locally splitting injection
on a neighborhood
of the base $B(C)\subset X_{{\mathbf F}_p}$.
Hence the assertion follows
from Proposition \ref{prsm}.

2.
For the $0$-section
$C=F^*T^*_XX|_{X_{{\mathbf F}_p}}$,
the base
$B(C)$ is $X_{{\mathbf F}_p}$
and the condition 
(1$''$) in the proof of
1 is satisfied.
Hence the assertion follows from 1.
\qed

}

\begin{pr}\label{prfCX}
Let $X,Y,V$
be regular noetherian schemes
satisfying the condition {\rm (F)}
and let
$$\begin{CD}
X@<h<<W\\
@VfVV@VV{f'}V\\
Y@<g<<V
\end{CD}$$
be a cartesian diagram
of morphisms of
finite type.
Let $C$
be a closed conical subset
of $FT^*X|_{X_{{\mathbf F}_p}}$.
Then the following conditions are
equivalent:

{\rm (1)}
$f$ is $C$-acyclic on
a neighborhood of the image
$h(W_{{\mathbf F}_p})$.

{\rm (2)}
There exists
a regular neighborhood $W_1\subset W$
of the inverse image of the base
$B(C)$ satisfying
the following conditions:
$f$ and $g$ are transversal
on $W_1$,
the restriction $h_1\colon W_1
\to X$ is $C$-transversal
and the restriction
$f'_1\colon W_1\to V$
is $h_1^\circ C$-acyclic.
\end{pr}

\proof{
First, 
we show that the both conditions
imply that $f$ is smooth on
a neighborhood of
the intersection $B(C)
\cap h(W_{{\mathbf F}_p})$.
For (1), this follows from
Lemma \ref{lmftr}.1.
For (2), similarly,
$f'_1$ is smooth on
a neighborhood of $B(h_1^\circ C)
=h_1^{-1}(B(C))$.
This implies that
$f$ is smooth on
a neighborhood of $
h(h_1^{-1}(B(C)))
=B(C)\cap h(W_{{\mathbf F}_p})$
since $f$ and $g$ are transversal
on $W_1$.

By replacing $X$ by
a neighborhood of
$B(C)\cap h(W_{{\mathbf F}_p})$
smooth over $Y$,
we may assume that 
$f$ is smooth.
Then, $W$ is regular
and $f$ and $g$ are transversal.
By Proposition \ref{prsm},
we have a commutative diagram
$$\begin{CD}
0\to& FT^*V|_{V_{{\mathbf F}_p}}
\times_{V_{{\mathbf F}_p}}
W_{{\mathbf F}_p}
@>>> FT^*W|_{W_{{\mathbf F}_p}}
@>>> F^*T^*W/V|_{W_{{\mathbf F}_p}}
&\to 0\\
&@AAA@AAA@AA{\cong}A&\\
0\to& FT^*Y|_{Y_{{\mathbf F}_p}}
\times_{Y_{{\mathbf F}_p}}
W_{{\mathbf F}_p}
@>>> FT^*X|_{X_{{\mathbf F}_p}}
\times_{X_{{\mathbf F}_p}}
W_{{\mathbf F}_p}
@>>> F^*T^*X/Y|_{X_{{\mathbf F}_p}}
\times_{X_{{\mathbf F}_p}}
W_{{\mathbf F}_p}
&\to0\\
\end{CD}$$
of exact sequences
of vector bundles on $W_{{\mathbf F}_p}$.
Let $C'\subset FT^*W|_{W_{{\mathbf F}_p}}$
be the image
of $h^*C=C\times_{X_{{\mathbf F}_p}}
W_{{\mathbf F}_p}
\subset FT^*X|_{X_{{\mathbf F}_p}}
\times_{X_{{\mathbf F}_p}}
W_{{\mathbf F}_p}$
and let
$A\subset FT^*Y|_{Y_{{\mathbf F}_p}}
\times_{Y_{{\mathbf F}_p}}
W_{{\mathbf F}_p}$
and $A'\subset FT^*V|_{V_{{\mathbf F}_p}}
\times_{V_{{\mathbf F}_p}}
W_{{\mathbf F}_p}$
be their inverse images.

Since the right vertical arrow is an
isomorphism,
the lower left arrow induces
an isomorphism
${\rm Ker}(FT^*Y|_{Y_{{\mathbf F}_p}}
\times_{Y_{{\mathbf F}_p}}
W_{{\mathbf F}_p}
\to 
FT^*V|_{V_{{\mathbf F}_p}}
\times_{V_{{\mathbf F}_p}}
W_{{\mathbf F}_p})
\to
{\rm Ker}(FT^*X|_{X_{{\mathbf F}_p}}
\times_{X_{{\mathbf F}_p}}
W_{{\mathbf F}_p}
\to FT^*W|_{W_{{\mathbf F}_p}})$.
Hence $A
\subset FT^*Y|_{Y_{{\mathbf F}_p}}
\times_{Y_{{\mathbf F}_p}}
W_{{\mathbf F}_p}$ is a subset
of the $0$-section
if and only if
$A'\subset
FT^*V|_{V_{{\mathbf F}_p}}
\times_{V_{{\mathbf F}_p}}
W_{{\mathbf F}_p}$ and
$h^*C\cap 
{\rm Ker}(FT^*X|_{X_{{\mathbf F}_p}}
\times_{X_{{\mathbf F}_p}}
W_{{\mathbf F}_p}
\to FT^*W|_{W_{{\mathbf F}_p}})$
are subsets
of the $0$-sections
and the assertion follows.
%
%
%
%
%
%
%
\qed

}

\begin{lm}\label{lmhf}
Let $f\colon X\to Y$
be a morphism of finite type of 
regular noetherian schemes
satisfying the condition {\rm (F)}.

{\rm 1.}
Let $C$ be a closed conical subset
of $FT^*X|_{X_{{\mathbf F}_p}}$
and 
assume that $f$ is proper
on the base $B(C)$.
Let $g\colon Y\to Z$
be a morphism of finite type of 
regular noetherian schemes
such that $Z_{{\mathbf F}_p}$
is of finite type over $k$.
Then the following conditions are
equivalent:

{\rm (1)}
$g$ is $f_\circ C$-acyclic.

{\rm (2)}
$gf$ is $C$-acyclic.

{\rm 2.}
Let $p\colon V\to X$
be a proper morphism of
regular schemes
and let $C=p_\circ F^*T^*_VV|_{V_{{\mathbf F}_p}}
\subset FT^*X|_{X_{{\mathbf F}_p}}$.
Then, the following conditions
are equivalent:

{\rm (1)}
$f$ is $C$-acyclic.

{\rm (2)}
The composition
$fp$ is smooth
on a neighborhood of $V_{{\mathbf F}_p}$.
\end{lm}

\proof{
1.
Let $x\in X_{{\mathbf F}_p}$ be a closed
point and $y\in Y_{{\mathbf F}_p}$ and $z\in Z_{{\mathbf F}_p}$
be the images.
Since the assertion is \'etale local,
we may also assume that
the morphisms 
$k(z)\to k(y)\to k(x)$
are isomorphisms.

Let $A_x$ be the inverse image of
$C_x$ by $F^*T^*_xX\gets F^*T^*_yY$.
Then, the inverse image $A'_x$
of $C_x$ by $F^*T^*_xX\gets F^*T^*_zZ$
equals the inverse image $A''_x$
of $A_x$
by $F^*T^*_yY\gets F^*T^*_zZ$.
Since the condition (1) 
(resp.\ (2)) is equivalent to
that $A'_x$ (resp.\ $A''_x$)
is a subset of the $0$-section
for any $x$,
the assertion follows.

2.
By 1.~applied to
$p_\circ F^*T^*_VV|_{V_{{\mathbf F}_p}}
=F^*T^*_VX|_{X_{{\mathbf F}_p}}$,
the condition (1)
is equivalent to that
the composition $fp$
is $F^*T^*_VV|_{V_{{\mathbf F}_p}}$-acyclic.
Hence the assertion follows from
Lemma \ref{lmftr}.2.
\qed

}

\begin{df}\label{dfhfC}
Let $X$
be a regular noetherian scheme
satisfying the condition {\rm (F)}
and let $C$ be a closed conical subset
of $FT^*X|_{X_{{\mathbf F}_p}}$.
We say that a pair
$(h,f)$ of morphisms
$h\colon W\to X$,
$f\colon W\to Y$
of finite type
of regular noetherian schemes
such that $Y_{{\mathbf F}_p}$
is of finite type over $k$
is $C$-acyclic
if the intersection of
$(C\times_{X_{{\mathbf F}_p}}
W_{{\mathbf F}_p})
\times_{W_{{\mathbf F}_p}} 
(FT^*Y|_{Y_{{\mathbf F}_p}}
\times_{Y_{{\mathbf F}_p}}
W_{{\mathbf F}_p})
\subset 
(FT^*X|_{X_{{\mathbf F}_p}}
\times_{X_{{\mathbf F}_p}}
W_{{\mathbf F}_p})
\times_{W_{{\mathbf F}_p}} 
(FT^*Y|_{Y_{{\mathbf F}_p}}
\times_{Y_{{\mathbf F}_p}}
W_{{\mathbf F}_p})$
with the kernel 
${\rm Ker}((h^*,f^*)\colon$
$(FT^*X|_{X_{{\mathbf F}_p}}
\times_{X_{{\mathbf F}_p}}
W_{{\mathbf F}_p})
\times_{W_{{\mathbf F}_p}} 
(FT^*Y|_{Y_{{\mathbf F}_p}}
\times_{Y_{{\mathbf F}_p}}
W_{{\mathbf F}_p})
\to
FT^*W|_{W_{{\mathbf F}_p}})$
is a subset of the $0$-section.
\end{df}

\begin{lm}\label{lmhfC}
Let $X$
be a regular noetherian scheme
satisfying the condition {\rm (F)}
and let $C$ be a closed conical subset
of $FT^*X|_{X_{{\mathbf F}_p}}$.

{\rm 1.}
Let $f\colon X\to Y$
be a morphism of finite type
of regular noetherian schemes
satisfying the condition {\rm (F)}.
Then, the following conditions are
equivalent:

{\rm (1)}
$f$ is $C$-acyclic.

{\rm (2)}
$(1_X,f)$ is $C$-acyclic.

{\rm 2.}
Let $h\colon W\to X$
and $f\colon W\to Y$
be morphisms of finite type
of regular noetherian schemes
satisfying the condition {\rm (F)}.
Then the following conditions are
equivalent:

{\rm (1)}
$(h,f)$ is $C$-acyclic.

{\rm (2)}
$h$ is $C$-transversal
and 
$f$ is $h^\circ C$-acyclic.
\end{lm}

\proof{1.
Identify the kernel of
$(1,f^*)\colon
FT^*X|_{X_{{\mathbf F}_p}}
\times_{X_{{\mathbf F}_p}} 
(FT^*Y|_{Y_{{\mathbf F}_p}}
\times_{Y_{{\mathbf F}_p}}
X_{{\mathbf F}_p})
\to
FT^*X|_{X_{{\mathbf F}_p}}$
with the image of
the injection
$(f^*,-1)\colon 
FT^*Y|_{Y_{{\mathbf F}_p}}
\times_{Y_{{\mathbf F}_p}}
X_{{\mathbf F}_p}
\to
FT^*X|_{X_{{\mathbf F}_p}}
\times_{X_{{\mathbf F}_p}}
(FT^*Y|_{Y_{{\mathbf F}_p}}
\times_{Y_{{\mathbf F}_p}}
X_{{\mathbf F}_p})$.
Then the inverse image
in $
FT^*Y|_{Y_{{\mathbf F}_p}}
\times_{Y_{{\mathbf F}_p}}
X_{{\mathbf F}_p}$ of
$C
\times_{X_{{\mathbf F}_p}}
(FT^*Y|_{Y_{{\mathbf F}_p}}
\times_{Y_{{\mathbf F}_p}}
X_{{\mathbf F}_p})
\subset
FT^*X|_{X_{{\mathbf F}_p}}
\times_{X_{{\mathbf F}_p}}
(FT^*Y|_{Y_{{\mathbf F}_p}}
\times_{Y_{{\mathbf F}_p}}
X_{{\mathbf F}_p})$
is the same as
the inverse image of
$C\subset T^*X$
and the assertion follows.

2.
Since ${\rm Ker}(h^*
FT^*X|_{X_{{\mathbf F}_p}}
\times_{X_{{\mathbf F}_p}}
W_{{\mathbf F}_p}
\to
FT^*W|_{W_{{\mathbf F}_p}})
\times 0
\subset
{\rm Ker}((h^*,f^*)\colon
(FT^*X|_{X_{{\mathbf F}_p}}
\times_{X_{{\mathbf F}_p}}
W_{{\mathbf F}_p})
\times_{W_{{\mathbf F}_p}} 
(FT^*Y|_{Y_{{\mathbf F}_p}}
\times_{Y_{{\mathbf F}_p}}
W_{{\mathbf F}_p})
\to
FT^*W|_{W_{{\mathbf F}_p}})$,
the $C$-acyclicity
of $(h,f)$ implies
the $C$-transversality of
$h$.
By 1.,
the $h^\circ C$-acyclicity
of $f$ is equivalent to the condition that
the intersection of
$h^\circ C\times_{W_{{\mathbf F}_p}} 
(FT^*Y|_{Y_{{\mathbf F}_p}}
\times_{Y_{{\mathbf F}_p}}
W_{{\mathbf F}_p})$
with 
${\rm Ker}(
FT^*W|_{W_{{\mathbf F}_p}}
\times_{W_{{\mathbf F}_p}}
(FT^*Y|_{Y_{{\mathbf F}_p}}
\times_{Y_{{\mathbf F}_p}}
W_{{\mathbf F}_p})
\to 
FT^*W|_{W_{{\mathbf F}_p}})$
is a subset of the
$0$-section.
This condition is equivalent to the
$C$-acyclicity
of $(h,f)$
since 
$h^\circ C\times_{W_{{\mathbf F}_p}} 
(FT^*Y|_{Y_{{\mathbf F}_p}}
\times_{Y_{{\mathbf F}_p}}
W_{{\mathbf F}_p})$
is the image of 
$h^*C\times_{W_{{\mathbf F}_p}} 
(FT^*Y|_{Y_{{\mathbf F}_p}}
\times_{Y_{{\mathbf F}_p}}
W_{{\mathbf F}_p})$.
\qed

}

\section{Micro-support}\label{sms}

We fix a perfect field $k$ 
of characteristic $p>0$
and a finite field 
$\Lambda$ of characteristic $\ell\neq p$.
We will assume
that a regular noetherian scheme $X$
over ${\mathbf Z}_{(p)}$
satisfies the condition {\rm (F)}
in Definition {\rm \ref{dfFTX}}.

\subsection{Micro-support}
\begin{df}\label{dfms}
Let $X$ be a regular noetherian scheme 
over ${\mathbf Z}_{(p)}$
satisfying the condition {\rm (F)}
in Definition {\rm \ref{dfFTX}} and
let $C$ be a closed conical subset
of the FW-cotangent bundle 
$FT^*X|_{X_{{\mathbf F}_p}}$.
Let ${\cal F}$
be a constructible complex
of $\Lambda$-modules
on $X$.
We say that ${\cal F}$
is micro-supported on $C$
if the following conditions
{\rm (1)} and {\rm (2)}
are satisfied:

{\rm (1)}
The intersection of
the support ${\rm supp}\, {\cal F}$
with the closed fiber $X_{{\mathbf F}_p}$ is 
a subset of the base $B(C)$.

{\rm (2)}
Every $C$-transversal separated morphism
$h\colon W\to X$ of finite type of
regular schemes
is ${\cal F}$-transversal
on a neighborhood of the closed
fiber $W_{{\mathbf F}_p}$.
\end{df}

This definition of micro-support
is related to
\cite[Proposition 8.13]{CC}
but is different from
\cite[1.3]{Be}.
We discuss this point in
Remark after Proposition \ref{prtrla}.
It is a property on a neighborhood
of $X_{{\mathbf F}_p}$.
If $X_{\mathbf Q}
=X\times_{{\rm Spec}\, {\mathbf Z}_{(p)}}
{\rm Spec}\, {\mathbf Q}$
is smooth over a field $K$
of characteristic $0$,
to cover $X_{\mathbf Q}$,
one can use the micro-support
of the restriction of ${\cal F}$ on $X_{\mathbf Q}$
defined as closed conical subset
of the cotangent bundle
$T^*X_{\mathbf Q}/K$.

\begin{lm}\label{lmTX}
Let $X$ be a regular noetherian scheme 
over ${\mathbf Z}_{(p)}$
satisfying the condition {\rm (F)}
and ${\cal F}$
be a constructible complex
of $\Lambda$-modules.

{\rm 1.}
${\cal F}$ is micro-supported
on $FT^*X|_{X_{{\mathbf F}_p}}$.

{\rm 2.}
If ${\cal F}$ is locally constant
on a neighborhood of
the closed fiber $X_{{\mathbf F}_p}$,
then 
${\cal F}$ is micro-supported
on the $0$-section $F^*T^*_XX|_{X_{{\mathbf F}_p}}$.

{\rm 3.}
Assume that $X$ is a
smooth scheme over $k$.
Let $C\subset T^*X$ be
a closed conical subset
and $F^*C\subset F^*T^*X
=FT^*X$ 
be the pull-back of $C$
Then, ${\cal F}$ is micro-supported
on $C$ in the sense of
{\rm (\cite[1.3]{Be}, \cite[Definition 4.1]{CC})}
if and only if 
${\cal F}$ is micro-supported
on $F^*C$.
\end{lm}

We show the converse of
2 in Corollary \ref{cortrla}.

\proof{
1.
Let $h\colon W\to X$
be a separated
morphism of finite type of regular schemes.
If $h$ is $FT^*X|_{X_{{\mathbf F}_p}}$-transversal,
then $h$ is smooth
on a neighborhood
of $W_{{\mathbf F}_p}$ by Lemma \ref{lmTXC}.1.
Hence $h$ is ${\cal F}$-transversal
on a neighborhood
of $W_{{\mathbf F}_p}$ by Lemma \ref{lmPoi}.1.

2.
Let $h\colon W\to X$
be a separated
morphism of finite type 
of regular schemes.
Then, since 
${\cal F}$ is locally constant
on a neighborhood of
the closed fiber $X_{{\mathbf F}_p}$,
$h$ is ${\cal F}$-transversal
on a neighborhood of $W_{{\mathbf F}_p}$
by Lemma \ref{lmPoi}.2.

3.
Let $h\colon W\to X$
be a separated
morphism of finite type of regular schemes.
Then, $h\colon W\to X$
is a separated morphism 
of smooth schemes
of finite type over $k$.
The morphism
$h\colon W\to X$ is $F^*C$-transversal
if and only if
$h\colon W\to X$ is $C$-transversal.
Hence
the equivalence follows
from \cite[Proposition 8.13]{CC}.
\qed

}

\begin{pr}\label{prmcf}
Let $X$ be a regular scheme 
over ${\mathbf Z}_{(p)}$
satisfying the condition {\rm (F)}
and let
${\cal F}$ be a constructible
complex of $\Lambda$-modules.
Let $C$ be a closed conical
subset of $FT^*X|_{X_{{\mathbf F}_p}}$
such that ${\cal F}$
is micro-supported on $C$.

{\rm 1.}
Let $h\colon W\to X$
be a separated morphism
of finite type of regular schemes.
If $h$ is $C$-transversal,
then $h$ is ${\cal F}$-transversal
on a neighborhood of $W_{{\mathbf F}_p}$
and $h^*{\cal F}$
is micro-supported on $h^\circ C$.

{\rm 2.}
Let $f\colon X\to Y$
be a separated
morphism of finite type 
proper on the base $B(C)$
of regular quasi-excellent
noetherian schemes
satisfying the condition {\rm (F)}.
Then $Rf_*{\cal F}$
is micro-supported on $f_\circ C$.
\end{pr}

\proof{
1.
Let $g\colon V\to W$
be an $h^\circ C$-transversal
separated morphism of finite type of
regular noetherian schemes.
Then, by Lemma \ref{lmhC},
$hg$ and $h$ are $C$-transversal.
Since ${\cal F}$ is
micro-supported on $C$,
$hg$ and $h$ are 
${\cal F}$-transversal
on neighborhoods of
$V_{{\mathbf F}_p}$ and of $W_{{\mathbf F}_p}$ respectively.
Hence by Proposition \ref{prhF}.1,
$g$ is $h^*{\cal F}$-transversal
on a neighborhood of
$V_{{\mathbf F}_p}$.

2.
Let $g\colon V\to Y$
be an $f_\circ C$-transversal 
separated morphism of finite type 
of regular noetherian schemes
and let
$$\begin{CD}
X@<h<< W\\
@VfVV@VV{f'}V\\
Y@<g<<V
\end{CD}$$
be a cartesian diagram.
Then, $f$ and $g$ are
transversal on a regular neighborhood
$W_1\subset W$
of the inverse image of
$B(C)$
and $h_1=h|_{W_1}\colon W_1\to X$
is $C$-transversal by Proposition \ref{prfC}.
Since ${\cal F}$
is micro-supported on $C$,
the restriction $h_1\colon W_1\to X$ is
${\cal F}$-transversal.

Since the intersection of
${\rm supp}\, {\cal F}$
with $X_{{\mathbf F}_p}$
is a subset of $B(C)$,
the intersection of
$A={\rm supp}\, h^*{\cal F}$
with $W_{{\mathbf F}_p}$ is a subset of $W_1$.
Since the closed
set $A\sm A\cap W_1$
does not intersect 
the closed fiber $W_{{\mathbf F}_p}$,
the complement $V_0=
V\sm f'(A\sm A\cap W_1)$
is an open neighborhood of $V_{{\mathbf F}_p}$.
By replacing $V$ by $V_0$,
we may assume $A={\rm supp}\, h^*{\cal F}\subset W_1$.
Then, $h\colon W\to X$ is
${\cal F}$-transversal.
Since $f$ and $g$ are transversal on $W_1$,
the base change morphism
$f'^*Rg^!\Lambda
\to Rh^!\Lambda$
is an isomorphism on $W_1$ by 
Lemma \ref{lmtrbc}.
Hence
$g$ is $Rf_*{\cal F}$-transversal
on $V$ by Corollary \ref{cortrC}.1.
\qed

}

We show that
being micro-supported is
an \'etale local property.

\begin{lm}\label{lmmset}
Let $X$ be a regular noetherian scheme 
over ${\mathbf Z}_{(p)}$
satisfying the condition {\rm (F)}
and let
${\cal F}$ be a constructible
complex of $\Lambda$-modules.

{\rm 1.}
Let $C$ be a closed conical
subset of $FT^*X|_{X_{{\mathbf F}_p}}$
and let
$(f_i\colon U_i\to X)_{i\in I}$
be an \'etale covering.
Then, the following conditions
are equivalent:

{\rm (1)}
${\cal F}$
is micro-supported on $C$.

{\rm (2)}
For every $i\in I$,
the pull-back
${\cal F}_i=f_i^*{\cal F}$
is micro-supported on 
$C_i=f_i^\circ C$.

{\rm 2.}
Let $f\colon U\to X$
be an \'etale morphism
of finite type
and let 
${\cal F}_U=f^*{\cal F}$ be the restriction.
Let
$C'\subset FT^*U$
be a closed conical subset
and let $C=f_*C'$ be
the closed conical subset defined
in Definition {\rm \ref{dfCet}}.
If ${\cal F}_U$ is micro-supported
on $C'$,
then ${\cal F}$ is micro-supported
on $C$.
\end{lm}

\proof{1.
The equivalence for the condition 
(1) in Definition \ref{dfms}
on the support is verified easily.
We show the equivalence
for the condition 
(2) in Definition \ref{dfms}
on transversality.

(1)$\Rightarrow$(2):
Let $i\in I$
and let $h\colon W\to U_i$ be a
separated morphism of finite type
of regular schemes.
If $h$ is $C_i$-tranversal,
then since $f_i$ is \'etale,
the composition
$f_i\circ h$ is $C$-tranversal
by Lemma \ref{lmhC}.
Hence
$f_i\circ h$ is ${\cal F}$-tranversal
on a neighborhood of the closed
fiber
and consequently
$h$ is ${\cal F}_i$-tranversal
on a neighborhood of the closed
fiber.

(2)$\Rightarrow$(1):
Let $h\colon W\to X$ be a
separated morphism of finite type
of regular schemes.
Assume $h$ is $C$-tranversal.
Then, for every $i\in I$,
the base change
$h_i\colon W_i=W\times_XU_i\to U_i$ 
is $C_i$-tranversal
by Lemma \ref{lmhC}.
Hence
$h_i$ is ${\cal F}_i$-tranversal
on a neighborhood of the closed
fiber
for every $i\in I$
and consequently
$h$ is ${\cal F}$-tranversal
on a neighborhood of the closed
fiber.

2.
The condition 
(1) in Definition \ref{dfms}
on the support is verified easily.
We show that the condition 
(2) in Definition \ref{dfms}
on transversality is satisfied.
Let $h\colon W\to X$
be a separated morphism
of finite type of regular
noetherian schemes.
Assume that $h\colon W\to X$
is $C$-transversal.
Then, by Lemma \ref{lmCet},
the morphism
$h\colon W\to X$
is smooth on a neighborhood $W_1$
of $h^{-1}(X_{{\mathbf F}_p}\sm 
f(U_{{\mathbf F}_p}))$.
Hence
$h|_{W_1}\colon W_1\to X$
is ${\cal F}$-transversal.
Further, by Lemma \ref{lmCet},
the base change
$W\times_XU\to U$ of $h$
is $C'$-transversal.
Since ${\cal F}_U$ is micro-supported
on $C'$,
the morphism
$W\times_XU\to U$
is ${\cal F}_U$-transversal
on a neighborhood of the closed
fiber.
Since $W_1$ and
$W\times_XU$
form an \'etale covering of $W$,
the morphism
$h$ is ${\cal F}$-transversal
on a neighborhood of the closed
fiber.
\qed

}

\medskip
In the rest of this subsection,
let ${\cal O}_K$ be
a discrete valuation ring
with residue field $k$
such that the fraction field
$K$ is of characteristic $0$.
Recall that ${\cal O}_K$ is excellent.
For a scheme $X$ over
${\cal O}_K$,
the closed fiber
$X_k=X\times_{{\rm Spec}\, 
{\mathcal O}_K}{\rm Spec}\, k$
has the same underlying set
as
$X_{{\mathbf F}_p}
=X\times_{{\rm Spec}\, 
{\mathbf Z}}{\rm Spec}\, 
{\mathbf F}_p$.

\begin{pr}\label{prtrla}
Let $h\colon W\to X$ and
$f\colon W\to Y$
be morphisms of regular schemes of finite
type over ${\cal O}_K$.
Let ${\cal F}$ be
a constructible complex
of $\Lambda$-modules
and $C$ be a closed
conical subset of $FT^*X|_{X_k}$.
Suppose that ${\cal F}$ is micro-supported
on $C$.
If the pair $(h,f)$ is
$C$-acyclic,
then 
$f\colon W\to Y$ is
$h^*{\cal F}$-acyclic along $W_k$.
\end{pr}

\proof{
By Lemma \ref{lmhfC}.2,
$h\colon W\to X$
is $C$-transversal
and $f\colon W\to Y$
is $h^\circ C$-acyclic.
Since 
${\cal F}$ is micro-supported on $C$,
the pull-back
$h^*{\cal F}$ is micro-supported
on $h^\circ C$ by Proposition \ref{prmcf}.1.
Hence by replacing $X$ by $W$,
we may assume $W=X$.

Since $f\colon X\to Y$ is
$C$-acyclic,
the morphism $f$ is
smooth on a neighborhood
of the intersection 
$B(C) \cap X_k
\supset
{\rm supp}\, {\cal F}
\cap X_k$ by Lemma \ref{lmftr}.1.
Hence, we may assume
$f\colon X\to Y$ is smooth.

Let $V\to Y$ be a separated morphism
of regular schemes
of finite type over ${\cal O}_K$.
Then the projection $p\colon U=V\times_YX\to X$
is $C$-transversal
by 
Proposition \ref{prfCX}. 
Hence $p$ is ${\cal F}$-transversal
on a neighborhood of $U_k$.
Thus by Proposition \ref{prla1},
$f$ is ${\cal F}$-acyclic
along $X_k$.
\qed

}

\medskip

\begin{rmk}\label{rmk}
{\rm 
The conclusion of Proposition
\ref{prtrla} is an analogue
of the original condition defining
the micro-support in \cite[1.3]{Be}.
In the geometric case,
this is shown to be equivalent
in \cite[Proposition 8.13]{CC}
to the condition analogous
to that in Definition \ref{dfms}.
However the following example
shows that the condition is
too weak in the setting of this article.

Let $X$ be a smooth scheme over
${\cal O}_K$ 
and $C=F^*T^*_{X_k}X|_{X_k}$
be the conormal bundle
of the closed fiber.
Let $(h,f)$ be a $C$-acyclic pair
of morphisms of regular schemes
of finite type over ${\cal O}_K$.
Then, since $h\colon W\to X$
is transversal to the immersion
$X_k\to X$,
the closed fiber $W_k$ is regular
and $W$ is smooth over ${\cal O}_K$
on a neighborhood of $W_k$
by Lemma \ref{lmhfC}.2
and Corollary \ref{corfC}.
Since $f\colon W\to Y$
is $F^*T^*_{W_k}W|_{W_k}$-acyclic,
further by Lemma \ref{lmhfC}.2
and Lemma \ref{lmhf}.2,
the morphism $W\to Y$ is smooth
on a neighborhood of $W_k$
and $W_k\to Y$ is also smooth.
This means that $W_k$ is empty.
Thus any ${\cal F}$ satisfies
the conclusion of
Proposition \ref{prtrla}.}
\end{rmk}

\begin{cor}\label{cortrla}
Let $X$ be a regular scheme of finite
type over ${\cal O}_K$
and ${\cal F}$
be a constructible complex
of $\Lambda$-modules.
Then, the following conditions
are equivalent.

{\rm (1)}
${\cal F}$ is locally constant
on a neighborhood of
the closed fiber $X_k$.

{\rm (2)}
${\cal F}$ is micro-supported
on the $0$-section $F^*T^*_XX|_{X_k}$.
\end{cor}

\proof{
(1)$\Rightarrow$(2)
is proved in Lemma \ref{lmTX}.2.

(2)$\Rightarrow$(1):
By Proposition \ref{prtrla}
applied to $(1_X,1_X)$,
the identity $1_X\colon X\to X$
is ${\cal F}$-acyclic
along $X_k$.
Hence ${\cal F}$ is locally constant
on a neighborhood of $X_k$
by Lemma \ref{lmla}.2.
\qed

}

%
%
%

\subsection{Singular support}

\begin{df}\label{dfss}
Let $X$ be a regular noetherian scheme 
over ${\mathbf Z}_{(p)}$
satisfying the condition {\rm (F)}
in Definition {\rm \ref{dfFTX}}
and ${\cal F}$
be a constructible complex
of $\Lambda$-modules
on $X$.
We say that a closed
conical subset $C
\subset FT^*X|_{X_{{\mathbf F}_p}}$ 
of the FW-cotangent bundle is the singular
support $SS{\cal F}$ of ${\cal F}$
if the following condition is satisfied:
For any closed conical subset
$C'\subset FT^*X|_{X_{{\mathbf F}_p}}$,
${\cal F}$ is micro-supported
on $C'$ if and only if $C\subset C'$.
\end{df}

If the singular support
$SS{\cal F}$ exists,
it is the intersection of
closed conical subsets $C
\subset FT^*X|_{X_{{\mathbf F}_p}}$ on
which
${\cal F}$ is micro-supported. 
The author does not know how
to prove the existence of
the singular support in general.
We compute 
the singular supports
in some cases.

\begin{lm}\label{lmsslc}
Let $X$ be a regular noetherian scheme 
over ${\mathbf Z}_{(p)}$
such that $X_{{\mathbf F}_p}$
is of finite type over $k$
and ${\cal F}$
be a constructible complex
of $\Lambda$-modules
on $X$.
We consider the following condition:

{\rm (DVR)}
$X$ is of finite type over ${\cal O}_K$
for a discrete valuation field $K$
of characteristic $0$
with residue field $k$.

{\rm 1.}
We consider
the following conditions:

{\rm (1)}
$SS{\cal F}$ is the $0$-section
$FT^*_XX|_{X_{{\mathbf F}_p}}$.

{\rm (2)}
${\cal F}$ is locally constant
on a neighborhood of $X_{{\mathbf F}_p}$
and $X_{{\mathbf F}_p}$
is a subset of the support of ${\cal F}$.

\noindent
We have
{\rm (2)}$\Rightarrow${\rm (1)}.
If {\rm (DVR)} is satisfied,
we have
{\rm (1)}$\Rightarrow${\rm (2)}.

{\rm 2.}
We consider
the following conditions:

{\rm (1)}
$SS{\cal F}=\varnothing$.

{\rm (2)}
${\cal F}=0$
on a neighborhood of $X_{{\mathbf F}_p}$.

\noindent
We have
{\rm (2)}$\Rightarrow${\rm (1)}.
If {\rm (DVR)} is satisfied,
we have
{\rm (1)}$\Rightarrow${\rm (2)}.

{\rm 3.}
Assume that 
$X={\rm Spec}\, {\cal O}_K$
for a discrete valuation ring
${\cal O}_K$ as in {\rm (DVR)}.
If ${\cal F}$
is not locally constant,
we have
$SS{\cal F}=FT^*X|_{X_k}$.
\end{lm}

Recall that a discrete valuation ring
${\cal O}_K$ as in {\rm (DVR)}
is excellent
by \cite[Scholie (7.8.3)]{EGA4}.

\proof{
1.
(1)$\Rightarrow$(2):
Since ${\cal F}$ is micro-supported
on the $0$-section
$F^*T^*_XX|_{X_{{\mathbf F}_p}}$,
by Corollary \ref{cortrla},
${\cal F}$ is locally constant
on a neighborhood of $X_{{\mathbf F}_p}$.
After replacing $X$ by
a neighborhood of $X_{{\mathbf F}_p}$,
we may assume that
${\cal F}$ is locally constant.
Then, the support
$Z={\rm supp}\, {\cal F}$
is an open and closed subset of $X$
and ${\cal F}$ is micro-supported
on the $0$-section $T^*_ZX|_{X_{{\mathbf F}_p}}$
on $Z$.
By the minimality of the singular
support, we have $T^*_ZX|_{X_{{\mathbf F}_p}}
=T^*_XX|_{X_{{\mathbf F}_p}}$
and $X_{{\mathbf F}_p}\subset Z$.

(2)$\Rightarrow$(1):
Since ${\cal F}$ is locally constant
on a neighborhood of $X_{{\mathbf F}_p}$,
by Lemma \ref{lmTX}.2,
${\cal F}$ is micro-supported
on the $0$-section
$F^*T^*_XX|_{X_{{\mathbf F}_p}}$.
Suppose ${\cal F}$ is micro-supported
on a closed conical subset
$C\subset FT^*X|_{X_{{\mathbf F}_p}}$.
Since $X_{{\mathbf F}_p}\subset {\rm supp}\, {\cal F}$,
we have
$X_{{\mathbf F}_p}\subset B(C)$.
This is equivalent to
$F^*T^*_XX|_{X_{{\mathbf F}_p}}\subset C$
and we obtain
$F^*T^*_XX|_{X_{{\mathbf F}_p}}=SS{\cal F}$.

2.
(1)$\Rightarrow$(2):
Since the intersection
${\rm supp}\, {\cal F}\cap
X_{{\mathbf F}_p}$ is a subset
of $SS{\cal F}=\varnothing$,
we have $X_{{\mathbf F}_p}\subset 
X\sm {\rm supp}\, {\cal F}$
and  the condition (2) holds.

(2)$\Rightarrow$(1):
Since every separated
morphism $h\colon W\to X$
is ${\cal F}$-transversal
on a neighborhood of $W_{{\mathbf F}_p}$
and since the intersection
${\rm supp}\, {\cal F}\cap
X_{{\mathbf F}_p}$ is empty,
${\cal F}$ is micro-supported
on $\varnothing$.

3.
By Lemma \ref{lmTX},
${\cal F}$ is micro-supported
on $FT^*X|_{X_k}$.
Suppose that ${\cal F}$ 
is micro-supported on 
a closed conical subset $C
\subset FT^*X|_{X_k}$.
Since $FT^*X|_{X_k}$
is a line bundle
over the point
${\rm Spec}\, k$,
$C$ is either
$\varnothing$,
the $0$-section
$F^*T^*_XX|_{X_k}$
or
$FT^*X|_{X_k}$ itself.
Since ${\cal F}$ is not locally
constant, 
by the contraposition of
Corollary \ref{cortrla} (2)$\Rightarrow$(1),
${\cal F}$ is not micro-supported
on the $0$-section
$F^*T^*_XX|_{X_k}$.
\qed

}

\medskip

We show that
being singular support is
a local property.

\begin{lm}\label{lmSSet}
Let $X$ be a regular noetherian scheme 
over ${\mathbf Z}_{(p)}$
satisfying the condition {\rm (F)}
and let
${\cal F}$ be a constructible
complex of $\Lambda$-modules.
Let $C$ be a closed conical
subset of $FT^*X|_{X_{{\mathbf F}_p}}$
and let
$(f_i\colon U_i\to X)_{i\in I}$
be an \'etale covering.
We consider the following conditions:

{\rm (1)}
$C=SS{\cal F}$.

{\rm (2)}
For every $i\in I$
and
the pull-backs
$C_i=f_i^\circ C$
and
${\cal F}_i=f_i^*{\cal F}$,
we have
$C_i=SS{\cal F}_i$.

\noindent
We have
{\rm (2)}$\Rightarrow${\rm (1)}.
If 
$(f_i\colon U_i\to X)_{i\in I}$
is a Zariski covering,
we have
{\rm (1)}$\Rightarrow${\rm (2)}
conversely.
\end{lm}

\proof{
Since the equivalence
of the condition 
to be micro-supported
is proved in Lemma \ref{lmmset},
we show the inplications
for the minimality.

(1)$\Rightarrow$(2):
Let $i\in I$
and let $C'\subset FT^*U_i$ be a
closed conical subset
on which ${\cal F}_i$
is micro-supported.
Then, by Lemma \ref{lmmset},
${\cal F}$ is micro-supported
on $f_{i*}C'$
in the notation of Definition \ref{dfCet}.
Hence, we have
$f_{i*}C'\supset SS{\cal F}=C$.
If $f_i\colon U_i\to X$
is an open immersion,
we have
$C'\supset f_i^\circ C=C_i$.

(2)$\Rightarrow$(1):
Let $C'\subset FT^*X$ be a
closed conical subset
on which ${\cal F}$
is micro-supported.
Then,
for every $i\in I$,
the pull-back
${\cal F}_i$ is micro-supported
on $f_i^\circ C'$
by Lemma \ref{lmmset}.
Hence, we have
$f_i^\circ C'\supset f_i^\circ C$.
Since $(U_i\to X)_{i\in I}$
is an \'etale covering,
we have $C'\supset C$.
\qed

}

\medskip

We compute the singular supports
of some sheaves
on regular schemes 
of finite type over a
discrete valuation ring.

\begin{lm}\label{lmXW}
Let $K$ be a discrete
valuation field of characteristic $0$
such that
the residue field $k$
is a perfect field of characteristic $p>0$.
Let $h\colon W\to X$ be a 
finite surjective morphism of
regular flat schemes
of finite type over
${\cal O}_K$
such that the morphism
$W_K\to X_K$
on the generic fiber is \'etale.
Assume that the reduced parts
$D=X_{k,{\rm red}}$
and
$E=W_{k,{\rm red}}$
of the closed fibers
are irreducible and are smooth 
of dimension $\geqq 1$
over the residue field $k$.

Assume that the following condition
is satisfied:

{\rm (1)}
The cokernel of the canonical morphism
$F\Omega^1_X
\otimes_{{\cal O}_X}
{\cal O}_E
\to
F\Omega^1_W
\otimes_{{\cal O}_W}
{\cal O}_E$
of locally free 
${\cal O}_E$-modules
is locally free of rank $1$.

{\rm 1.}
The direct image
$C= \pi_\circ FT^*_WW|_E
\subset FT^*X|_D$
of the $0$-section
is the image of the sub line bundle
${\rm Ker}(FT^*X|_D\times_DE
\to FT^*W|_E)$
of
$FT^*X|_D\times_DE$.

%
%
%
%
%
%
%
{\rm 2.}
Further assume 
that the following condition is satisfied:

{\rm (2)}
The finite morphism
$E\to D$ is purely inseparable
of degree $\geqq 1$.

\noindent
Then, for each closed point $x\in D$
and for the point $w\in E$
above $x$,
there exists
a regular subscheme
$Z\subset W$ 
of codimension $1$
containing $w$ and
flat over ${\cal O}_K$
satisfying the following conditions:

The composition
$Z\to W\to X$ is unramified.
The pull-back $C\times_{X_{{\mathbf F}_p}}w
\subset FT^*X\times_{X_{{\mathbf F}_p}}
w$ 
of the fiber at $x$
equals the fiber of the
kernel of the surjection
$FT^*X\times_{X_{{\mathbf F}_p}}Z
\to FT^*Z$.
\end{lm}

\proof{
1. 
Since the ${\cal O}_E$-linear morphism
$F\Omega^1_X
\otimes_{{\cal O}_X}
{\cal O}_E
\to
F\Omega^1_W
\otimes_{{\cal O}_W}
{\cal O}_E$
of locally free 
${\cal O}_E$-modules
of the same rank has
the cokernel of rank 1,
the kernel is also locally free of
rank 1.
Hence the assertion follows.

2.
Let $n=\dim {\cal O}_{X,x}$.
Since $E\to D$ is assumed
to be purely inseparable,
the residue field
$k(w)$ is a purely inseparable
extension of a perfect field
$k(x)$ and hence 
the morphism $k(x)\to k(w)$ is an isomorphism.
By the assumption on the
rank of the cokernel
and by Proposition \ref{prdx},
the $k(x)$-linear mapping
${\mathfrak m}_x/
{\mathfrak m}_x^2
\to
{\mathfrak m}_w/
{\mathfrak m}_w^2$
induced by
${\cal O}_{X,x}\to
{\cal O}_{W,w}$ is of rank $n-1$.

Take an element of
${\mathfrak m}_w/
{\mathfrak m}_w^2$
not contained in the image
of ${\mathfrak m}_x/
{\mathfrak m}_x^2$
and take its lifting
$f\in {\mathfrak m}_w$
not divisible by
a prime element $t$ defining 
the divisor $E\subset W$.
Then,
a regular closed subscheme $Z$
of codimension $1$ 
of a neighborhood 
of $w$ is defined 
by $f$.
Let $z$ denote $w\in W$
regarded as a point of $Z$.
Since $f$ is not divisible by $t$,
we may assume that $Z$ is flat
over ${\cal O}_K$.

Since $\bar f\in
{\mathfrak m}_w/
{\mathfrak m}_w^2$ is 
not contained in the image
of ${\mathfrak m}_x/
{\mathfrak m}_x^2$,
the induced morphism
${\mathfrak m}_x/
{\mathfrak m}_x^2
\to
{\mathfrak m}_w/((f)+
{\mathfrak m}_w^2)
=
{\mathfrak m}_z/
{\mathfrak m}_z^2$
is a surjection.
Hence further shrinking 
$Z$ if necessary,
we may assume that
$Z\to X$ is unramified.
Since the kernel of the surjection
${\mathfrak m}_x/
{\mathfrak m}_x^2
\to
{\mathfrak m}_z/
{\mathfrak m}_z^2$
equals the kernel of
${\mathfrak m}_x/
{\mathfrak m}_x^2
\to
{\mathfrak m}_w/
{\mathfrak m}_w^2$,
the last condition 
on the fibers is satisfied.
\qed

}

\medskip
We show that some concrete
examples of Kummer coverings
satisfy the assumptions
in Lemma \ref{lmXW}.
Let $K$ be a discrete valuation
field as in Lemma \ref{lmXW}
containing a primitive
$p$-th root of 1.
Let $X$ be a regular flat
scheme of finite type
over ${\cal O}_K$
and assume that the reduced part
$D=X_{k,{\rm red}}$
is smooth over the residue field $k$.
Let $L$ be the local field 
at the generic point of $D$
and let $e={\rm ord}_Lp\geqq p-1$
be the absolute ramification index.

\begin{lm}\label{lmKum}
Let $\pi \in \Gamma(X,{\cal O}_X)$
be a uniformizer of the divisor $D
=X_{k,{\rm red}}\subset X$
and let $u \in \Gamma(X,{\cal O}_X^\times)$
be a unit.
Let $1\leqq n< \dfrac{pe}{p-1}$
be an integer
congruent to $0$ or $1$
modulo $p$
and set $n=pm$ or $n=pm+1$
respectively.
In the case $n=pm$,
assume that $du$ defines locally
a part of a basis of $\Omega^1_D$.
Define a Kummer covering
$V\to U=X_K$
by $v^p=1+u\pi^n$.

{\rm 1.}
The normalization $\pi\colon
W\to X$
in $V$ is regular.
The reduced closed fiber
$E=W_{k,{\rm red}}$
is smooth over $k$
and the finite morphism
$E\to D$ is purely inseparable.

{\rm 2.}
The cokernel
${\rm Coker}(F\Omega^1_X
\otimes_{{\cal O}_X}
{\cal O}_E
\to
F\Omega^1_W
\otimes_{{\cal O}_W}
{\cal O}_E)$
is an invertible ${\cal O}_E$-module.

{\rm 3.}
Assume $n=pm$.
If $e=m+1$,
let $\pi'$ denote the uniformizer
$p/\pi^m$.
Then, the kernel of the canonical morphism
$FT^*X|_D\times_DE\to
FT^*W|_E$
is a line bundle spanned
by 
$$\begin{cases}
w(u)-u\cdot w(\pi')
&
\text{ if $p=2$ and $e=m+1$},
\\
w(u)&
\text{ otherwise}.
\end{cases}$$
\end{lm}

\proof{
{\rm 1.}
Since the assertion is local,
we may assume that
$X={\rm Spec}\, A$ is affine.
We show that the normalization
$B$ of $A$ is generated by 
$t=(v-1)/\pi^m$.
By the assumption
$n<\dfrac{ep}{p-1}$,
we have
$e+m>pm$ and
the polynomial
$(1+\pi^mT)^p-1
\in A[T]$
is divisible by $\pi^{pm}$.
Define a monic polynomial
$F\in A[T]$
by $1+\pi^{pm}F=(1+\pi^mT)^p$.
Since
$F\equiv T^p\bmod 
\pi A$ and since $u$ is a unit,
in the case $n=pm+1$,
the equation
$F=\pi u$ is an Eisenstein equation.
In the case $n=pm$,
the reduction of the equation
$F=u$
modulo $\pi A$
gives $T^p=u$.
In this case $du$ is a part of
a basis of $\Omega^1_D$
by the assumption.
Hence 
by setting $v=1+\pi^mt$
where $t\in B$ denotes the class of
$T$,
we obtain
$B=A[T]/(F-u\pi)$ 
or $B=A[T]/(F-u)$ 
respectively.

The reduced part $E$
is defined by $t$ or $\pi$
according to $n=pm+1$ or $n=pm$
respectively.
Hence $E$ is smooth over $k$
and the finite morphism
$E\to D$ is purely inseparable
of degree 1 or $p$
respectively.

2.
By Corollary \ref{corXZ},
we have a commutative diagram
$$\begin{CD}
0@>>>
F^*N_{D/X}
\otimes_{{\cal O}_D}
{\cal O}_E
@>>>
F\Omega^1_X
\otimes_{{\cal O}_X}
{\cal O}_E
@>>>
F^*\Omega^1_D
\otimes_{{\cal O}_D}
{\cal O}_E
@>>>
0\\
@.@VVV@VVV@VVV@.\\
0@>>>
F^*N_{E/W}
@>>>
F\Omega^1_W
\otimes_{{\cal O}_W}
{\cal O}_E
@>>>
F^*\Omega^1_E
@>>>
0
\end{CD}$$
of exact sequences of
locally free ${\cal O}_E$-modules.
In the case $n=pm+1$,
the right vertical arrow
is an isomorphism
since $E\to D$ is an isomorphism.
Further
the left vertical arrow is $0$
since the ramification index is $p$.
In the case $n=pm$,
the left vertical arrow
is an isomorphism
since the ramification index is $1$.
Further
the cokernel of the right vertical arrow is 
locally free of rank 1
since $E\to D$ is a purely inseparable
covering defined by $T^p=u$
and $du$ is a part of a basis
of $\Omega^1_D$.
Hence the assertion follows.

3.
We compute
the polynomial $F
\bmod  \pi^2$.
Recall that we have
$e+m>pm$.
Since $e$ is divisible by
$p-1$, the equality
$e+m=pm+1$ holds
if and only if
$p=2$ and $e=m+1$.
Hence 
the coefficients of $T^i$ 
for $i=1,\ldots, p-1$ in 
the polynomial $F$
are divisible by $\pi^2$
except $F=T^2+
2/\pi^m\cdot T$
in the exceptional case.

Thus, except the exceptional case,
we have a congruence
$F\equiv T^p
\bmod \pi^2$
and hence the kernel is
spanned by 
$w(u)$.
In the exceptional case,
we have
$t^2+\pi't=u$ for
$\pi'=2/\pi^m$.
Hence $w(u)$
is sent to $t^2\cdot w(\pi')
=u\cdot w(\pi')$.
\qed

}

\begin{pr}\label{prKum}
Let $K$ be a discrete
valuation field of characteristic $0$
such that
the residue field $k$
is a perfect field of characteristic $p>0$.
Let $X$ be a regular flat scheme
of finite type over
${\cal O}_K$
such that the reduced part
$D=X_{k,{\rm red}}$
is irreducible and is smooth
over the residue field $k$.

Let ${\cal F}_U$ be a locally constant
constructible sheaf of
$\Lambda$-modules on
the generic fiber $U=X_K$
and let ${\cal F}=j_!{\cal F}_U$
be the $0$-extension
for the open immersion
$j\colon U\to X$.
Let $V\to U$ be a finite
\'etale Galois covering 
of Galois group $G$ such that
the pull-back ${\cal F}_V$
is constant 
and let $\pi\colon W\to X$ be
the normalization 
of $X$ in $V$.

Assume that $W$ is regular
and that
the reduced part
$E=W_{k,{\rm red}}$
is also irreducible and smooth
over the residue field $k$.
Assume that the order of $G$
is invertible in $\Lambda$
and that ${\cal F}_U$ corresponds
to a non-trivial
irreducible representation $M$ of $G$.

{\rm 1.}
The canonical morphism
${\cal F}=j_!{\cal F}_U
\to Rj_*{\cal F}_U$
is an isomorphism.

{\rm 2.}
Assume that conditions
{\rm (1)} and {\rm (2)} 
in Lemma {\rm \ref{lmXW}}
are satisfied.
Then, the singular support $SS{\cal F}$
equals the direct image
$C=\pi_\circ FT^*_WW|_{W_k}$
of the $0$-section.
\end{pr}

\proof{
1.
By the assumption 
that the order of $G$
is invertible in $\Lambda$
and that $M$ is an irreducible
representation,
the locally constant sheaf
${\cal F}_U$
is isomorphic to a direct summand
of $\pi_{K*}\Lambda$
where $\pi_K\colon V=W_K\to U=X_K$
is the restriction of $\pi$.

Let $j_W\colon W_K\to W$
be the open immersion
of the generic fiber.
Since $W$ is regular
and the reduced part
of the closed fiber
$W_k$ is a regular divisor,
we have isomorphisms
$\Lambda\to j_{W*}\Lambda$,
$\Lambda_E(-1)\to R^1 j_{W*}\Lambda$
and $R^qj_{W*}\Lambda=0$
for $q\neq 0,1$
by the absolute purity 
\cite[{\sc Th\'eor\`eme 3.1.1}]{purete}.
Similarly,
we have isomorphisms
$\Lambda\to j_*\Lambda$,
$\Lambda_D(-1)\to R^1 j_*\Lambda$
and $R^qj_*\Lambda=0$
for $q\neq 0,1$.
Since $E\to D$ induces a homeomorphism
on the \'etale site by the assumption,
the canonical morphism
$\Lambda_D\to \pi_*\Lambda_E$
is an isomorphism.
Hence, for the cokernel
${\cal G}={\rm Coker}
(\Lambda_X\to \pi_*\Lambda_W)$,
the canonical morphisms
$j_!j^*{\cal G}\to
{\cal G}\to Rj_*j^*{\cal G}$
are isomorphisms.

Since
$M$ is a non-trivial irreducible
representation
of a semi-simple algebra
$\Lambda[G]$,
the corresponding sheaf
${\cal F}$ is a direct summand
of $j^*{\cal G}$.
Hence the canonical morphism
${\cal F}=j_!{\cal F}_U
\to Rj_*{\cal F}_U$
is an isomorphism.

2.
Since ${\cal F}$ is a direct summand
of $\pi_*\Lambda_W=\Lambda_X
\oplus {\cal G}$,
by Proposition \ref{prmcf}.2,
the constructible sheaf
${\cal F}$ is micro-supported on
$C=\pi_\circ FT^*_WW|_{W_k}$.

Suppose ${\cal F}$
is micro-supported on 
a closed conical subset $C'$.
It suffices to prove $C\subset C'$.
Let $x\in X_{{\mathbf F}_p}$
be a closed point,
let $h\colon Z\to X$
be an unramified morphism
as in Lemma \ref{lmXW}
and let $z\in Z$ be
the unique point above $x$.
Since $Z\to X$ factors through
$Z\to W$,
the restriction 
${\cal F}_{Z\cap U}$
is constant.
Hence the morphism
$h$ is not
${\cal F}$-transversal
by the contraposition
of Proposition \ref{prhF}.2
(1)$\Rightarrow$(2)
and 1.
Since ${\cal F}$
is micro-supported on $C'$,
the morphism
$h$ is not $C'$-transversal,
on any open neighborhood of $z\in Z$.

The kernel $L={\rm Ker}
(FT^*X\times_{X_{{\mathbf F}_p}}
Z_{{\mathbf F}_p}\to FT^*Z)$
is a line bundle
on $Z_{{\mathbf F}_p}$.
The intersection $C'_1=h^*C'
\cap L
\subset FT^*X\times_{X_{{\mathbf F}_p}}
Z_{{\mathbf F}_p}$
is a closed conical subset
of $L$.
Let $Z_1
=\{y\in Z_{{\mathbf F}_p}\mid
C'_{1,y}=L_y\}$
be the image by the projection 
of the complement
$C'_1\sm (C'_1\cap Z_{{\mathbf F}_p})$
of the $0$-section.
Since $C'_1\subset L$ is a closed
conical subset,
the image
$Z_1\subset
Z_{{\mathbf F}_p}$ is a closed subset.
Since the restriction 
$Z\sm Z_1\to X$ of $h$
is $C'$-transversal,
the complement
$Z\sm Z_1$ is not
an open neighborhood of $z$.
Namely,
we have $z\in Z_1$
and hence 
$C'_{1,z}=L_z$ 
is a subset of $C'_z$.

Since 
$L_z=C_z=C_x\times_xz$
by the condition on $Z$, we get
$C_x\subset C'_x$
for each closed point $x\in X_k$.
Thus we have $C\subset C'$
as required.
\qed

}


\begin{thebibliography}{99}

\bibitem{cst}
M.\ Artin,
{\em Faisceaux constructibles
Cohomologie d'une courbe
alg\'ebrique},
SGA 4 Expos\'e IX,
Th\'eorie des Topos et Cohomologie \'Etale des Sch\'emas, 
Lecture Notes in Mathematics Volume 305, 1973, pp 1-42.


\bibitem{Artin}
-----,
{\em Morphismes acycliques},
SGA 4 Expos\'e XV,
Th\'eorie des Topos et Cohomologie \'Etale des Sch\'emas, 
Lecture Notes in Mathematics Volume 305, 1973, pp 168-205.

\bibitem{smbc}
-----,
{\em Th\'eor\`eme de changement de base
par morphisme lisse},
SGA 4 Expos\'e XVI,
Th\'eorie des Topos et Cohomologie \'Etale des Sch\'emas, 
Lecture Notes in Mathematics Volume 305, 1973, pp 206-249.

\bibitem{Be}
A.\ Beilinson,
{\em Constructible sheaves are holonomic},
Selecta Mathematica,
volume 22, 1797-1819 (2016).

\bibitem{dJ}
A.\ J.\ de Jong,
{\em Smoothness, semi-stability and alterations},
Publ.\ Math.\ IHES, 83,
(1996),  51--93.


\bibitem{DP}
P.\ Deligne,
{\em La formule de dualit\'e globale},
Th\'eorie des topos et
cohomologie \'etale des sch\'emas,
SGA 4 Expos\'e XVIII,
Springer Lecture Notes in Math.\ 305 (1972), 480-587.

\bibitem{Rapport}
------, 
{\em Rapport sur la formule des traces},
Cohomologie \'etale SGA 4$\frac12$,
Springer Lecture Notes in Math.\ 569 (1977)
76--109.



\bibitem{TF}
------, 
{\em Th\'eor\`emes de finitude en cohomologie $\ell$-adique},
Cohomologie \'etale SGA 4$\frac12$,
Springer Lecture Notes in Math.\ 569 (1977)
233--251.


\bibitem{DKRZ}
T.\ Dupuy, E.\ Katz, J.\ Rabinoff, D.\ Zureick-Brown,
{\em Total $p$-differential on schemes over 
${\mathbf Z}/p^2$}, 
Journal of Algebra 524, 110-123 (2019).



\bibitem{GR}
O.\ Gabber, L.\ Ramero,
{\em Foundations for almost ring theory -- Release 7.5},
{\tt https://}
{\tt arxiv.org/abs/math/0409584}.

\bibitem{EGA4}
A.\ Grothendieck,
{\sc \'El\'ements de g\'eom\'etrie alg\'ebrique IV}, 
\'Etude locale des sch\'emas et
des morphismes de sch\'emas, 
Publ.\ Math.\ IHES 20, 24, 28, 32 (1964-67).



\bibitem{Ill}
L.\ Illusie,
{\em Complexe cotangent et d\'eformations I},
Springer Lecture Notes in Math., 239,
Springer-Verlag, Berlin, Heidelberg, New York 1971.

\bibitem{app}
-----, 
{\em Appendice \`a Th\'eor\`emes de finitude en cohomologie $\ell$-adique},
Cohomologie \'etale SGA 4$\frac12$,
Springer Lecture Notes in Math.\ 569 (1977)
252--261.


\bibitem{KSc}
M.\ Kashiwara, P.\ Schapira,
{\sc Sheaves on manifolds},
Springer-Verlag, Grundlehren der Math.\ Wissenschaften 292,
(1990).


\bibitem{fini}
F.~Orgogozo, 
{\em Le th\'eor\`eme de finitude}, XIII,
Travaux de Gabber
sur l'uniformisation locale et la cohomologie \'etale des sch\'emas quasi-excellents,
Ast\'erisque 363-364 (2014)
261-276.

\bibitem{purete}
J.~Riou, 
{\em Classes de Chern, morphismes de Gysin, puret\'e absolue}, XVI,
Travaux de Gabber
sur l'uniformisation locale et la cohomologie \'etale des sch\'emas quasi-excellents,
Ast\'erisque 363-364 (2014)
301-350.

\bibitem{CC}
T.\ Saito,
{\em The characteristic cycle and the singular support of a constructible sheaf},
Invent.\ Math.\ (2017) 207, 597--695.

\bibitem{gr}
------,
{\em Graded quotients of ramification groups
of local fields with imperfect residue fields},
{\tt arXiv:2004.03768}

\bibitem{FW}
------,
{\em Frobenius-Witt differentials and regularity},
{\tt arXiv:2008.04728}

\end{thebibliography}
\end{document}